\documentclass[11pt]{report}

\usepackage{graphicx}
\usepackage{doublespace}

\begin{document}
\pagenumbering{roman}
\begin{center}
{\Large \bf Green's Function of 3-D Helmholtz Equation for Turbulent Medium: Application to Optics}
\centerline{By}
\centerline{\large  PENG LI}
\centerline{\large B.A. (Fudan University) 2001}
\centerline{\large  THESIS}
\centerline{\large Submitted in partial satisfaction of the requirements for the
degree of}
\centerline{\large  MASTER OF ARTS}
\centerline{\large in}
\centerline{\large APPLIED MATHEMATICS}
\centerline{\large in the}
\centerline{\large  OFFICE OF GRADUATE STUDIES}
\centerline{\large of the}
\centerline{\large  UNIVERSITY OF CALIFORNIA}
\centerline{\large  DAVIS}
\end{center}
{\large Approved:}
\begin{center}
\centerline{Prof. Albert Fannjiang\ \underbar{\hskip 2.5in}}
\centerline{Prof. Naoki Saito\ \ \ \ \ \ \ \ \underbar{\hskip 2.5in}}
\centerline{Prof. Hong Xiao\ \ \ \ \ \ \ \ \ \underbar{\hskip 2.5in}}
\centerline{\large Committee in Charge}
\centerline{\large 2004}
\end{center}
\newpage
\large

\newpage
\large
{\Large \bf ACKNOWLEDGMENTS} \\
\addcontentsline{toc}{chapter}{\numberline{}ACKNOWLEDGMENTS}

I would like to start by thanking my advisor Professor Albert C. Fannjiang of the Mathematics Department at the University of California, Davis. He has been very supportive and allowed me the freedom to explore many diverse areas of study. He has also inspired many ideas in the topic of waves in inhomogeneous media and provided invaluable advices on my research.\\
Next I would like to thank Professor Naoki Saiko and Professor Hong Xiao of the mathematics department at UC-Davis. Both of them are experts in computational mathematics and have helped me greatly in the area of the apodization problem and its corresponding numerical analysis. I also appreciate their being my thesis committee members.\\
Last I would like to thank some of my professors and fellow students both from UC-Davis and from Fudan University where I completed my B.A. in applied math. Professor Zhaojun Bai who recommended me to the Graduate Group in Applied Mathematics at UC-Davis. Professor Alexander Soshnikov, Thomas Strohmer, Roger Wets, Steve Shkoller, Janko Gravner for their in-class instruction as well as many fruitful discussions. Another person I am greatly thankful to is our graduate coordinator Celia Davis, who has consistently supported me during my study at UC-Davis. Finally I would like to thank my parents in China. Without their support I would not have made any progress in any aspect. 

\newpage 
\begin{center}
\underline{\bf \Large Abstract} 
\end{center}
\addcontentsline{toc}{chapter}{\numberline{}Abstract}

\medskip

\large
The fundamental problem of optical wave propagation is the determination of the field at an observation point, given a disturbance specified over some finite aperture. In both vacuum and inhomogeneous media, the solution of this problem is given approximately by the superposition integral, which is a mathematical expression of the extended Huygens-Fresnel principle.\\

In doing so, it is important to find the atmospheric impulse response (Green's function). Within a limited but useful region of validity, a satisfactory optical propagation theory for the earth's clear turbulent atmosphere could be developed by using Rytov's method to approximate the Helmholtz equation. In particular, we deal with two optical problems which are the time reversal and apodization problems. The background and consequences of these results for optical communication through the atmosphere are briefly discussed.

\begin{spacing}{2}

\tableofcontents

\newpage
\pagestyle{myheadings} 
\pagenumbering{arabic}

\markright{  \rm \normalsize CHAPTER 1. \hspace{0.5cm}
  Introduction}
\large 

\chapter{Introduction}
\thispagestyle{myheadings}

In contemporary radiophysics, atmospheric optics and wireless communication, one often studies the propagation of electromagnetic waves in the atmosphere. In doing so, it is increasingly important to take into account the turbulent state of the atmosphere, a state which produces fluctuations in the refractive index of the air.\\

Determining the turbulence effect on wave propagation requires one to solve Maxwell's equations with a refractive index that is stochastic. To solve for effects of turbulence in the low-power regime, one starts with Maxwell's equations
\begin{eqnarray}
\nabla \times \vec{\varepsilon} (\vec{r}, t)= - \frac{\partial \vec{B} (\vec{r}, t)}{\partial t},\\
\nabla \cdot \vec{\varepsilon} (\vec{r}, t)= \frac{q(\vec{r}, t)}{\epsilon (\vec{r}, t)},\\
\nabla \times \vec{B} (\vec{r}, t) = \mu_0 \vec{\rho} (\vec{r}, t) + \mu_0 \epsilon (\vec{r}, t) \frac{\partial \vec{\varepsilon} (\vec{r}, t)}{\partial t}, and\\
\nabla \cdot \vec{B} (\vec{r}, t) =0,
\end{eqnarray}
where $\vec{\varepsilon} (\vec{r}, t)$ is the electric field, $\vec{B} (\vec{r}, t)$ is the magnetic field, $\epsilon (\vec{r}, t)$ is the permittivity, $\mu_0$ is the magnetic permeability, $\vec{\rho} (\vec{r}, t)$ is the current density, and $q(\vec{r}, t)$ is the charge density. Assume there is a time-harmonic variation of the electric field
\begin{equation}
\vec{\varepsilon} (\vec{r},t)=exp(-i\omega t) \vec{E}(\vec{r})
\end{equation}
where $\vec{E}(\vec{r})$ can also vary with time but the time scale of variation is much longer than the harmonic variation. The air density affects the refractive index as
\begin{equation}
\epsilon (\vec{r}) = \epsilon_0 n^2 (\vec{r})
\end{equation}
Then the wave equation describing the electric field derived from the above relations is
\begin{equation}
\Delta \vec{E}(\vec{r})+k_0^2 n^2(r)\vec{E}(\vec{r})-2i \frac{k_0}{c} \frac{\partial [n^2(\vec{r}) \vec{E}(\vec{r})]}{\partial t}+2 \nabla \{\vec{E}(\vec{r}) \cdot \nabla ln[n(\vec{r})]\}=0
\end{equation}
where $k_0=\frac{\omega}{c} = \frac{2\pi}{\lambda}$ is the free space wavenumber with $c$ being the speed of light in vacuum.\\

By reasonable simplifications (details in Chapter 4), one finally need to study the Helmholtz equation (reduced wave equation) for the turbulent medium
\begin{equation}
\Delta \vec{E}(\vec{r}) +k_0^2n^2(\vec{r}) \vec{E}(\vec{r})=0
\end{equation}
with appropriate boundary conditions describing the evolution of the harmonic amplitude for the components of the electric field.\\

Within a limited but useful region of validity, a satisfactory optical propagation theory has been developed by using Rytov's method to approximate the Helmholtz equation. If it is assumed that the magnitude of the air-density inhomogeneity is small, the refractive index is then
\begin{equation}
n(\vec{r}) = 1+\delta n_1(\vec{r})
\end{equation}
where $\delta n_1(\vec{r})<<1$.\\

By using the Rytov approximation, we can obtain not only the solution, but also the Green's function of the Helmholtz equation for the turbulent medium with the aid of the Green's function for free-space propagation.\\
It can also be shown that the Green's function in turbulence obtained by Rytov's method satisfies the same reciprocity condition as the free-space Green's function.\\ 

A better understanding of how randomness and turbulence affect the propagating wave fields can lead to better solutions to many problems in the optical society. Two examples that we are particularly interested when understanding wave turbulence interaction is important are:\\
I.) Time reversal problem.\\
II.) Apodization problem.\\

In both these two examples we are concerned with the relationship between the initial distribution of light over the exit pupil of an optical system and the amplitude distribution in the object plane. This is based on the extended Huygens-Fresnel Principle in random medium developed by H.T.Yura [12, 13]. The superposition integral requires an explicit form of the atmospheric impulse response (Green's function), which we have obtained by using Rytov's method. Further analysis is given for each problem.\\

The structure of this thesis is as the following:\\
We start Chapter 2 by constructing the Green's function for the free-space 3-D Helmholtz equation (\ref{eq:hel}) with the Sommerfeld radiation condition (\ref{eq:rad}). With a point source at a point $\vec{\xi}$, the Green's function is represented by (\ref{eq:qwe}).\\

In Chapter 3 we will give a brief introduction to the refractive index variation. We will also explain how the refractive index variation is determined by the temperature fluctuations and present some experimental results.\\

In Chapter 4 and 5, we first introduce the Rytov approximation method, and then derive the solution and Green's function for the Helmholtz equation in turbulent medium step by step. At the end of Chapter 5, the reciprocity property of the Green's function for general Helmholtz equation will be proved.\\

Further simplification by parabolic approximation is conducted in Chapter 6. When solving for the solution of Helmholtz equation, we specialize the initial condition into two cases: plane wave and beam wave. An explicit form of the Green's function will also be derived by parabolic approximation.\\

In Chapter 7 and 8, we will apply the results obtained in previous work to two important applications: time reversal problem and apodization problem.\\

A brief discussion is given in Chapter 9.

\markright{  \rm \normalsize CHAPTER 2. \hspace{0.5cm}
  Green's Function for Free-Space 3-D Helmholtz Equation}
\large 

\chapter{Green's Function for Free-Space 3-D Helmholtz Equation}
\thispagestyle{myheadings}

We consider the 3-D Wave Equation\\
\begin{equation} \label{eq:wav}
-(U_{tt}-c^2\Delta U)=q(\vec{r},t)
\end{equation}
where $q(\vec{r},t)$ is a source.
If $q(\vec{r},t)=q(\vec{r})exp(-i\omega t)$ represents a source oscillating with a single frequency $\omega$, then the entire motion reduces to a wave motion with same frequency $\omega$ after an initial transient period, so that we can write
\begin{equation}
U(\vec{r},t)=U(\vec{r})exp(-i\omega t)
\end{equation}
Thus (\ref{eq:wav}) reduces to the 3-D Helmholtz equation
\begin{equation} \label{eq:hel}
-(\Delta+k^2)U(\vec{r})=f(\vec{r})
\end{equation}
where $k=\frac{\omega}{c}, f(\vec{r})=c^{-2}q(\vec{r})$\\

The function $U(\vec{r})$ satisfies this equation in some domain $D\subset R$ with boundary $\partial D$, and it also satisfies some prescribed boundary conditions. We also assume that $U(\vec{r})$ satisfies the Sommerfeld radiation condition
\begin{equation} \label{eq:rad}
\lim_{r \to \infty}r(U_r-ikU)=0
\end{equation}
which simply states that the solution behaves like outgoing waves generated by the source.\\

We construct a Green's function $G(\vec{r},\vec{\xi})$ for (\ref{eq:hel}), so that $G(\vec{r},\vec{\xi})$ satisfies the equation
\begin{equation} \label{eq:gre}
-(\Delta + k^2)G=\delta(x)\delta(y)\delta(z)
\end{equation}\\

Using the spherical polar coordinates, the 3-D Laplacian can be expressed in terms of radial coordinate r,
\begin{equation} 
\Delta G=G_{rr}+\frac{2}{r} G_r
\end{equation}
so that (\ref{eq:gre}) assumes the form
\begin{equation}
-[\frac{1}{r^2} \frac{\partial}{\partial r} (r^2 \frac{\partial G}{\partial r})+k^2 G]=\delta (r), \ \ \ \ 0<r<\infty
\end{equation}
with the radiation condition (\ref{eq:rad}).\\

For r$>$0, the function G satisfies the homogeneous equation
\begin{equation}
\frac{1}{r^2} \frac{\partial}{\partial r} (r^2 \frac{\partial G}{\partial r})+k^2 G=0
\end{equation}
Or, equivalently,
\begin{equation}
\frac{\partial^2}{\partial r^2} (rG)+k^2(rG)=0
\end{equation}

This equation admits a solution of the form
\begin{equation}
rG(r)=Ae^{ikr}+Be^{-ikr}
\end{equation}
or
\begin{equation} \label{eq:sln}
G(r)=A\frac{e^{ikr}}{r}+B\frac{e^{-ikr}}{r}
\end{equation}
where A and B are arbitrary constants. In order to satisfy the radiation condition, we need to set B=0, and hence, solution (\ref{eq:sln}) becomes
\begin{equation}
G(r)=A\frac{e^{ikr}}{r}
\end{equation}

To determine A, we use the spherical surface $S_\epsilon$ of radius $\epsilon$ and divergence theorem, so that
\begin{equation}
\lim_{\epsilon \to 0} \int_{S_\epsilon} \frac{\partial G}{\partial r} dS=\lim_{\epsilon \to 0} \int_{S_\epsilon} \frac{A}{r} e^{ikr}(ik-\frac{1}{r}) dS=-1
\end{equation}
from which we find $A=\frac{1}{4\pi}$ as $\epsilon \to 0$. Consequently, the Green's function takes the form
\begin{equation}
G(r)=\frac{e^{ikr}}{4\pi r}
\end{equation}

Physically, this represents outgoing spherical waves radiating away from the source at the origin. With a point source at a point $\vec{\xi}$, the Green's function is represented by
\begin{equation} \label{eq:qwe}
G(\vec{r},\vec{\xi})=\frac{exp\{ik|\vec{r}-\vec{\xi}|\}}{4\pi |\vec{r}-\vec{\xi}|}
\end{equation}
where $\vec{r}$ and $\vec{\xi}$ are position vectors in $R^3$.\\

Finally, when k=0, this result reduces exactly to the Green's function for the 3-D Poisson equation.

\newpage
\pagestyle{myheadings} 
\markright{  \rm \normalsize CHAPTER 3. \hspace{0.5cm} 
  Introduction to Refractive Index Variations}

\chapter{Introduction to Refractive Index Variations}
\thispagestyle{myheadings}

We are all familiar with the twinkling of stars at night and the shimmering of distant objects on a hot day. These visual effects are caused by small refractive index inhomogeneities in the atmosphere, which in turn are produced by atmospheric turbulence. Here we refer the word \textsl{turbulence} to the density fluctuations arising from atmospheric temperature fluctuations. When a beam of light passes through the air above a road heated by the sun, the randomly fluctuating air temperature produces small refractive index inhomogeneities that affect the beam propagation. \\

Consider for example an initially well-defined phase front propagating through a region of atmospheric turbulence. Because of random fluctuations in phase velocity the initially well defined phase front will become distorted. This alters and re-directs the flow of energy in the beam. As the distorted phase front progresses, random changes in beam direction and intensity fluctuations occur. The beam is also found to spread in size beyond the dimensions predicted by diffraction theory in free space.\\

\begin{figure}[h]
\begin{center}
\includegraphics[width=12cm]{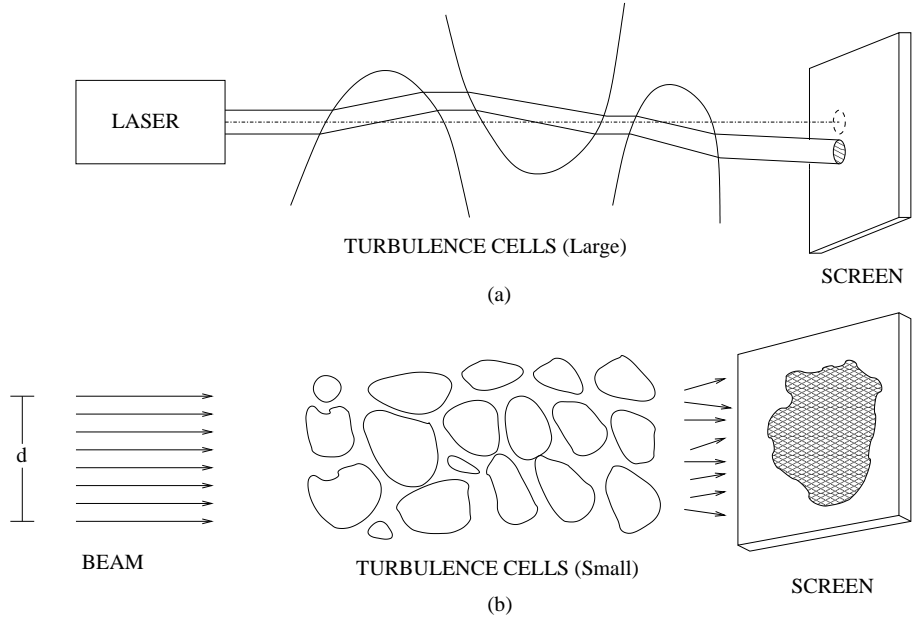}
\caption{(a) Laser beam is deflected by turbulence cells that are larger than the beam diameter, and (b) Laser beam is broken up by turbulence cells that are smaller than the beam diameter. } \label{fig:tur1}
\end{center}
\end{figure}

The reason why all this happens, as we have stated, is atmospheric turbulence that arises when air parcels of different temperatures are mixed by wind and convection. The individual air parcels, or turbulence cells, break up into smaller cells and eventually lose their identity. In the meantime, however, the mixing produces fluctuations in the density and therefore in the refractive index of air. A cartoon of laser beam propagation in turbulence is given in Figure (\ref{fig:tur1})\\

To describe these random processes, we need to find a way to define the fluctuations that are characteristic of turbulence. This is done by expressing each quantity as the sum of an average and a fluctuation term. We can write the random variation of the refractive index as
\begin{equation} \label{eq:or}
n(\vec{r})=n_0+\delta n(\vec{r}), \ \ \ with\ \ \ \ \ n_0=< n(\vec{r}) >
\end{equation}
Here the sharp brackets indicate the ensemble average
\begin{equation}
< n(\vec{r}) > = \frac{1}{\tau} \int_0^\tau n(\vec{r}, t) dt
\end{equation}
where $\tau$ is a time that is large compared to the lifetime of a fluctuation. A repetition of the averaging operation has no effect on $<n(\vec{r})>$, and it follows that $<[<n(\vec{r})>]> = <n(\vec{r})>$ and $<\delta n(\vec{r})>=0$.\\

In later chapters when we consider the Helmholtz equation in turbulence (details are given in Chapter 4)
\begin{equation}
\Delta u(\vec{r}) + k_0^2 n^2(\vec{r}) u(\vec{r})=0
\end{equation}
We can write
\begin{equation} \label{eq:ref}
n(\vec{r}) = 1+\delta n_1(\vec{r})
\end{equation}
since we can assume $< n(\vec{r}) >=1$ by absorbing a constant term into the wavenumber $k_0$, and write
\begin{equation}
\delta n(\vec{r}) = \delta \cdot n_1(\vec{r})
\end{equation}
where $\delta \cdot n_1(\vec{r})<<1$, and $\delta$ was inserted into this expression to show smallness.\\

Now let's return to the original expression (\ref{eq:or}). The refractive index fluctuation $\delta n(\vec{r})$ is proportional to the density fluctuation $\delta N(\vec{r})$, which in turn is governed by temperature fluctuations. The dependence of $\delta n(\vec{r})$ on $\delta N(\vec{r})$ follows from the Lorenz-Lorentz law
\begin{equation} \label{eq:l5.3}
\frac{n^2-1}{n^2+2} \cdot \frac{1}{N} = C_1
\end{equation}
where $C_1$ is a constant (for a given frequency of the electromagnetic wave) and $N$ is the particle density. For air, with $n\simeq 1$, equation (\ref{eq:l5.3}) may be approximated by
\begin{equation} \label{eq:l5.4}
(n-1)=\frac{3}{2} N C_1
\end{equation}

Differentiation of the above equation yields the relationship between $\delta n$ and $\delta N$:
\begin{equation}
\delta n=\frac{3}{2} C_1 \delta N
\end{equation}

If we assume that air obeys the perfect gas law, then $N=\frac{p}{kT}$ and it follows from (\ref{eq:l5.4}) that
\begin{equation}
n-1=\frac{3 C_1 p}{2 k T}
\end{equation}
where p is pressure expressed in millibars and T is temperature measured in Kelvins. For the optical region of the spectrum, the constant factor $\frac{3C_1}{2k}$ is given approximately by $79\times 10^{-6}$, thus
\begin{equation} \label{eq:l5.7}
(n-1)=79\times 10^{-6} \frac{p}{T}
\end{equation}

We see that the refractive index depends on both pressure and temperature. Since pressure variations are relatively small, fluctuations of the refraction index are primarily due to temperature fluctuations. Because their short lifetime and small-scale sizes may be treated as an adiabatic process, hence
\begin{equation} \label{eq:l5.8}
p T^{\frac{\gamma}{(1-\gamma)}}= constant
\end{equation}
where $\gamma \approx 1.4$ for air.\\

Differentiation of equation (\ref{eq:l5.8}) leads to
\begin{equation}
\frac{\delta p}{p} = (\frac{\gamma}{\gamma-1}) \frac{\delta T}{T}
\end{equation}

Now with the aid of (\ref{eq:l5.7}) we can eliminate $\delta p$ and find that the refraction index fluctuations are proportional to the temperature fluctuations, that is
\begin{equation}
\delta n = \frac{79 \times 10^{-6}}{\gamma-1} \frac{p}{T^2} \delta T
\end{equation}

The single most important parameter appearing in almost all equations that describe beam disturbances caused by turbulence is the refractive index structure coefficient $C_n$. It is governed by the pressure and temperature difference at two points separated by the distance $r$(measured in centimeters) and is given by
\begin{equation} \label{eq:l5.11}
C_n = [ 79 \times 10^{-6} \frac{p}{T^2} ] C_T
\end{equation}
where the temperature structure parameter is
\begin{equation}
C_T = \sqrt{<(T_1 - T_2)^2>} \ r^{-\frac{1}{3}}
\end{equation}
The temperatures $T$, $T_1$, and $T_2$ are all in Kelvins, and $p$ is the atmospheric pressure in millibars. For strong turbulence, $C_n=5 \times 10^{-7}$; for intermediate turbulence, $C_n=4 \times 10^{-8}$; and for weak turbulence, $C_n=8 \times 10^{-9}$. These are typical values. The structure parameter usually appears in equations in the form $C_n^2$, which varies from about $10^{-17} m^{-\frac{2}{3}}$ for extremely weak turbulence to $10^{-13} m^{-\frac{2}{3}}$ or more when the turbulence is strong. This latter value is usually observed near the ground in direct sunlight. Measurements of $C_n$ with temperature sensors 1.6m above the ground have shown that the minimum value of $C_n$ occurs about one to two hours before sunrise and after sunset. The peak values have been recorded around noon on sunny days. As might be expected form equation (\ref{eq:l5.11}), the structure parameter $C_n$ decreases with altitude. This is shown in Table 3.1.

\begin{center}
\begin{tabular}{ccc}
\multicolumn{2}{c}{Table 3.1   Typical values of $C_n$ as a function of height. 
 }\\
Height(Km) & $C_n (m^{-1/3}) \times 10^8$ \\ 0.001 & 30 \\ 0.003 & 20 \\0.01 & 15 \\ 0.03 & 10 \\ 0.1 & 6 \\ 0.3 & 4 \\ 1.0 & 1 \\ 3.0 & 1
\end{tabular}
\end{center}

\begin{figure}[h]
\begin{center}
\includegraphics[width=10cm]{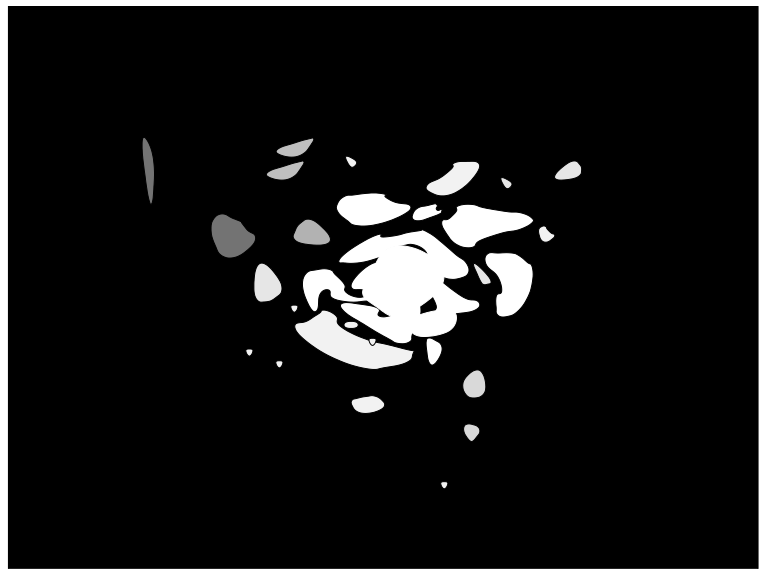}
\caption{Instantaneous intensity distribution of a laser beam at 15.5 km. Area of screen covered by photo is about 1.20 m by 1.55 m. } \label{fig:tur2}
\end{center}
\end{figure}

The movement of small index-of-refraction inhomogeneities through the path of beam causes random deflection and interference between different portions of the wavefront, which can lead to an internal breaking up of the beam spot into smaller \textsl{hot spots}. Figure (\ref{fig:tur2}) shows a typical instantaneous intensity distribution of a He-Cd laser beam at 15.5 km. The initial intensity distribution of this beam at the laser was Gaussian. The photo illustrates clearly the complexity of the transmitted pattern that consists of amorphous areas of random sizes and shapes. Other patterns may have a reticulated appearance with sharp lines bounding large polygonal areas. The dominant size of the hot spots in the pattern is given approximately by $\sqrt{L\lambda}$, where $L$ is the distance from the laser to the observation screen. These bright patches of about 1-cm diameter are typical for $L\approx 1km$. For space-to-earth paths, $L$ would be the distance between the observation screen and the turbulence layer in the atmosphere.\\

\begin{figure}
\begin{center}
\includegraphics[width=10cm]{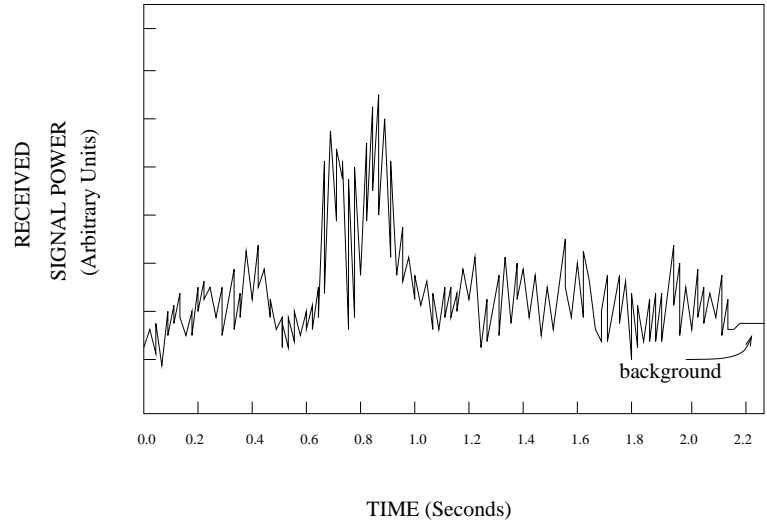}
\caption{Typical received signal fluctuations for a small detector 145 km from the transmitter. (From A. L. Buck, Applied Optics 6, 703, April 1967)} \label{fig:tur3}
\end{center}
\end{figure}

Because of the constantly changing pattern, a small detector placed in the beam will measure intensity fluctuations or scintillation. Typical power fluctuations recorded by a detector whose aperture size is considerably less than the beam diameter are shown in Figure (\ref{fig:tur3}). 


The temporal frequency of the intensity fluctuations recorded at a fixed point within the beam usually varies between 1 and 100 Hz.\\

In determining the turbulence effect on wave propagation, it was Rytov who proposed an approximation that included diffraction effects. This resulted in a wide range of validity and many problems have been solved as interest in wave propagation in turbulence expanded because of new technologies and areas of application. These areas include sending laser beams efficiently through the atmosphere, remoting sensing of the atmosphere on earth and correcting atmospheric distortion to allow better resolution in optical images.\\
When propagating laser beams parallel to the ground, one found that scintillation increased with increasing turbulence as the Rytov theory predicted until a certain level was reached at which the measured scintillation saturated. Fortunately, even when the scintillation is saturated Rytov's method typically gives the correct answer for phase disturbances, thus allowing one to treat many problems of practical interest. It has generally been found that the Rytov approximation gives a good approximation if
\begin{equation}
k_0^{\frac{7}{6}} \int_0^L C_n^2(z) z^{\frac{5}{6}} dz < 1
\end{equation}
where $k_0$ is the free-space wavenumber, $L$ is the distance from the laser to the observation screen and $C_n$ is the refractive index structure coefficient.

\newpage
\pagestyle{myheadings} 
\markright{  \rm \normalsize CHAPTER 4. \hspace{0.5cm} 
  Rytov Approximation I - Solution of Helmholtz Equation}

\chapter{Rytov Approximation I - Solution of Helmholtz Equation}
\thispagestyle{myheadings}

Determining the turbulence effect on wave propagation requires one to solve Maxwell's equations with a refractive index that is stochastic.\\
Assume there is a time-harmonic variation of the electric field
\begin{equation}
\vec{\varepsilon} (\vec{r},t)=exp(-i\omega t) \vec{E}(\vec{r})
\end{equation}
where $\vec{E}(\vec{r})$ can also vary with time but the time scale of variation is much longer than the harmonic variation.\\
The wave equation describing the electric field derived from the Maxwell's equation is then
\begin{equation} \label{eq:max}
\Delta \vec{E}(\vec{r})+k_0^2 n^2(r)\vec{E}(\vec{r})-2i \frac{k_0}{c} \frac{\partial [n^2(\vec{r}) \vec{E}(\vec{r})]}{\partial t}+2 \nabla \{\vec{E}(\vec{r}) \cdot \nabla ln[n(\vec{r})]\}=0
\end{equation}
where $k_0=\frac{\omega}{c} = \frac{2\pi}{\lambda}$ is the free space wavenumber with $c$ being the speed of light in vacuum, and $n(\vec{r})$ is the refractive index.\\

If the rate of change of $\vec{E}(\vec{r})$ with time is much less than the sinusoidal variation $exp(-i\omega t)$, then the third term is negligible compared to the first two. For optical frequencies this condition is satisfied.\\

In addition, if the propagation wavelength $\lambda$ is much less than the inner scale, we can ignore depolarization effects that is due to the last term in (\ref{eq:max}). This corresponds to scatterers being large relative to the wave length and the direction of the scattered wave very close to the direction of the original scattered wave which is the case for the regime that we will consider.
For visible wavelengths, typical inner-scale sizes of a millimeter satisfy this condition, in which case one obtains a scalar equation for each of the electric field components separately. Therefore in this regime the Helmholtz equation is 
\begin{equation}
\Delta \vec{E}(\vec{r}) +k_0^2n^2(\vec{r}) \vec{E}(\vec{r})=0
\end{equation}
with appropriate boundary conditions describing the evolution of the harmonic amplitude for the components of the electric field.\\

The equation for one component is
\begin{equation} \label{eq:1cm}
\Delta u(\vec{r}) + k_0^2 n^2(\vec{r}) u(\vec{r})=0
\end{equation}

In the Rytov method the solution is expressed as
\begin{equation}
u(\vec{r})=exp[\Phi(\vec{r})]
\end{equation}

This leads to the nonlinear Riccati equation
\begin{equation} \label{eq:Ric}
\Delta \Phi (r)+|\nabla \Phi (r)|^2 = - k_0^2 n^2(r)
\end{equation}

If it is assumed that the magnitude of the air-density inhomogeneity is small, the refractive index is then
\begin{equation} \label{eq:reff}
n(\vec{r}) = 1+\delta n_1(\vec{r})
\end{equation}
where $\delta n_1(\vec{r})<<1$, and $\delta$ was inserted into this expression to show smallness.\\

By perturbation method, we assume that the solution can be written as a power series in $\delta$, i.e.
\begin{equation} \label{eq:ser}
\Phi(\vec{r},\delta)=\Phi_0 (\vec{r}) + \delta \Phi_1(\vec{r})+ \delta^2 \Phi_2(\vec{r})+......
\end{equation}

By inserting (\ref{eq:reff}) and (\ref{eq:ser}) into (\ref{eq:Ric}), and by separating the equations based on the power of $\delta$, one obtains a system of equations 
\begin{equation} \label{eq:1st}
\delta^0: \Delta \Phi_0(\vec{r}) +\nabla \Phi_0 (\vec{r})\cdot \nabla \Phi_0 (\vec{r})=-k_0^2
\end{equation}
\begin{equation} \label{eq:2nd}
\delta^1: \Delta \Phi_1(\vec{r}) +2\nabla \Phi_0 (\vec{r})\cdot \nabla \Phi_1 (\vec{r})=-2k_0^2 n_1 (\vec{r})
\end{equation}
\begin{equation}
\delta^2: \Delta \Phi_2(\vec{r}) +2\nabla \Phi_0 (\vec{r})\cdot \nabla \Phi_2 (\vec{r})=-k_0^2 n_1 (\vec{r})-\nabla \Phi_1 (\vec{r})\cdot \nabla \Phi_1 (\vec{r})
\end{equation}
............
\begin{equation}
\delta^m: \Delta \Phi_m(\vec{r}) +2\nabla \Phi_0 (\vec{r})\cdot \nabla \Phi_m (\vec{r})=-\sum_{p=1}^{m-1}\nabla \Phi_p (\vec{r})\cdot \nabla \Phi_{m-p} (\vec{r})
\end{equation}

For these equations to be valid, one requires
\begin{equation}
|\nabla \Phi_{n+1}(\vec{r})|<<|\nabla \Phi_n(\vec{r})|
\end{equation}

We retain only the first two terms, and express the lowest-order term in the form
\begin{equation}
u_0(\vec{r})=exp[\Phi_0(\vec{r})]
\end{equation}
then equation (\ref{eq:1st}) is equivalent to the free-space wave equation
\begin{equation} \label{eq:fre}
\Delta u_0 (\vec{r}) + k_0^2 u_0(\vec{r})=0
\end{equation}

Now we set 
\begin{equation}
\Phi_1(\vec{r}) = \frac{W_1 (\vec{r})}{u_0(\vec{r})}
\end{equation}

By inserting this into equation (\ref{eq:2nd}) and by applying (\ref{eq:fre}), we need to obtain the solution of
\begin{equation}
\Delta W_1(\vec{r}) + k_0^2 W_1(\vec{r})=-2k_0^2 n_1(\vec{r}) u_0(\vec{r})
\end{equation}

The solution to this linear differential equation is then obtained with the Green's function for free-space propagation as 
\begin{equation}
\Phi_1(\vec{r})=\frac{2k_0^2}{u_0(\vec{r})}\int dV' u_0(\vec{r}') n_1(\vec{r}') G(|\vec{r}-\vec{r}'|)
\end{equation}
where the integration is over the source volume denoted by primed coordinates.\\
We have obtained the Green's function in Chapter 1, which is
\begin{equation}
G(r)=\frac{exp(ik_0r)}{4\pi r}
\end{equation}
Thus the perturbed field is 
\begin{equation} \label{eq:fin}
\Phi_1(\vec{r})=\frac{k_0^2}{2\pi u_0(\vec{r})}\int dV' u_0(\vec{r}') n_1(\vec{r}') \frac{exp(ik_0|\vec{r}-\vec{r}'|)}{|\vec{r}-\vec{r}'|}
\end{equation}

By (\ref{eq:ser}), the solution of (\ref{eq:1cm}) to the 2nd order is then
\begin{equation}
u(\vec{r})=u_0(\vec{r})exp[\delta \Phi_1(\vec{r})]
\end{equation}
where $\Phi_1(\vec{r})$ is given by (\ref{eq:fin}).

\newpage
\pagestyle{myheadings} 
\markright{  \rm \normalsize CHAPTER 5. \hspace{0.5cm} 
  Rytov Approximation II - Green's Function of Helmholtz Equation}

\chapter{Rytov Approximation II - Green's Function of Helmholtz Equation}
\thispagestyle{myheadings}

\section{Green's Function of Helmholtz Equation by Rytov Approximation}

Following the previous chapter, the Helmholtz equation in turbulence for one component is
\begin{equation} \label{eq:1cm1}
\Delta u(\vec{r}) + k_0^2 n^2(\vec{r}) u(\vec{r})=0
\end{equation}

We also assume that the magnitude of the inhomogeneity is small, so that the refractive index is 
\begin{equation} \label{eq:ref1}
n(\vec{r}) = 1+\delta n_1(\vec{r})
\end{equation}
where $\delta n_1(\vec{r})<<1$, and $\delta$ was inserted into this expression to show smallness.\\

Then the Green's function $G(\vec{r},\vec{\xi})$ for (\ref{eq:1cm1}) satisfies the following equation
\begin{equation}
\Delta G(\vec{r},\vec{\xi}) + k_0^2 n^2(\vec{r}) G(\vec{r},\vec{\xi})=-\delta (|\vec{r}-\vec{\xi}|)
\end{equation}

In the Rytov method, the solution is expressed as
\begin{equation}
G(\vec{r},\vec{\xi})=exp(\Phi (\vec{r},\vec{\xi}))
\end{equation}
which leads to the following equation
\begin{equation}
e^{\Phi (\vec{r},\vec{\xi})} [ \Delta\Phi (\vec{r},\vec{\xi})+\nabla \Phi (\vec{r},\vec{\xi}) \cdot \nabla \Phi (\vec{r},\vec{\xi})] + k_0^2 n^2 (\vec{r}) e^{\Phi (\vec{r},\vec{\xi})}= -\delta (|\vec{r}-\vec{\xi}|)
\end{equation}
or equivalently
\begin{equation} \label{eq:3.5}
\Delta\Phi (\vec{r},\vec{\xi})+\nabla \Phi (\vec{r},\vec{\xi}) \cdot \nabla \Phi (\vec{r},\vec{\xi}) = -k_0^2 n^2 (\vec{r}) - e^{- \Phi (\vec{r},\vec{\xi})}\delta (|\vec{r}-\vec{\xi}|)
\end{equation}

As in chapter 2, we assume that the solution of (\ref{eq:3.5}) can be written as a power series in $\delta$,
\begin{equation} \label{eq:3.6}
\Phi (\vec{r},\vec{\xi}, \delta)=\Phi_0 (\vec{r},\vec{\xi}) + \delta \Phi_1 (\vec{r},\vec{\xi})+ \delta^2 \Phi_2 (\vec{r},\vec{\xi}) + \dots
\end{equation}

Notice that the refractive index is given by (\ref{eq:ref1}), then the right hand side of equation (\ref{eq:3.5}) could be written as
\begin{equation}
RHS=-k_0^2 (1+2\delta n_1(\vec{r}) +\delta^2 n_1^2(\vec{r}))-e^{- \Phi_0 (\vec{r},\vec{\xi})} (1-\delta \Phi_1 (\vec{r},\vec{\xi})+ O (\delta^2)) \delta (|\vec{r}-\vec{\xi}|)
\end{equation}
where we have approximated $\Phi (\vec{r},\vec{\xi})$ by the first two terms in (\ref{eq:3.6}), and applied Tayler's expansion on $e^{-\delta \Phi_1 (\vec{r},\vec{\xi})}$.

By inserting (\ref{eq:3.6}) into (\ref{eq:3.5}) and by separating the equations based on the power of $\delta$, we obtained the equation for $\delta^0$ as
\begin{equation} \label{eq:3.8}
\Delta\Phi_0 (\vec{r},\vec{\xi}) + \nabla \Phi_0 (\vec{r},\vec{\xi}) \cdot \nabla \Phi_0 (\vec{r},\vec{\xi}) = 
-k_0^2 - e^{- \Phi_0 (\vec{r},\vec{\xi})} \delta (|\vec{r}-\vec{\xi}|)
\end{equation}

Now we set
\begin{equation}
G_0(\vec{r},\vec{\xi})=exp(\Phi_0 (\vec{r},\vec{\xi}))
\end{equation}
then equation (\ref{eq:3.8}) is equivalent to 
\begin{equation} \label{eq:3.9}
\Delta G_0 (\vec{r},\vec{\xi}) + k_0^2 G_0 (\vec{r},\vec{\xi}) = -\delta (|\vec{r}-\vec{\xi}|)
\end{equation}
which is the equation to solve for the Green's function for the free-space Helmholtz equation.\\

It has been showed in chapter1 that
\begin{equation} \label{eq:3.10}
G_0 (\vec{r},\vec{\xi})= \frac{e^{ik_0 |\vec{r}-\vec{\xi}|}}{4\pi |\vec{r}-\vec{\xi}|}
\end{equation}

The equation for $\delta^1$ is 
\begin{equation} \label{eq:3.11}
\Delta\Phi_1 (\vec{r},\vec{\xi}) + 2 \nabla \Phi_0 (\vec{r},\vec{\xi}) \cdot \nabla \Phi_1 (\vec{r},\vec{\xi}) = -2k_0^2 n_1(\vec{r}) + e^{- \Phi_0 (\vec{r},\vec{\xi})} \Phi_1 (\vec{r},\vec{\xi}) \delta (|\vec{r}-\vec{\xi}|)
\end{equation}

If we set
\begin{equation}
W_1 (\vec{r},\vec{\xi}) = \Phi_1 (\vec{r},\vec{\xi}) G_0 (\vec{r},\vec{\xi})
\end{equation}

then
\begin{equation} \label{eq:3.13}
\Delta W_1 (\vec{r},\vec{\xi})=G_0 (\vec{r},\vec{\xi}) \Delta \Phi_1 (\vec{r},\vec{\xi}) + 2 G_0 (\vec{r},\vec{\xi}) \nabla \Phi_0 (\vec{r},\vec{\xi}) \cdot \nabla \Phi_1 (\vec{r},\vec{\xi}) + \Phi_1 (\vec{r},\vec{\xi})\Delta G_0 (\vec{r},\vec{\xi})
\end{equation}

Now by combining (\ref{eq:3.9}), (\ref{eq:3.11}) and (\ref{eq:3.13}), we obtain
\begin{eqnarray}
\Delta W_1 (\vec{r},\vec{\xi}) + k_0^2 W_1 (\vec{r},\vec{\xi}) 
\nonumber &=& G_0 (\vec{r},\vec{\xi}) (\Delta \Phi_1 (\vec{r},\vec{\xi}) + 2 \nabla \Phi_0 (\vec{r},\vec{\xi}) \cdot \nabla \Phi_1 (\vec{r},\vec{\xi}))\\ \nonumber &\ & + \Phi_1 (\vec{r},\vec{\xi}) (\Delta G_0 (\vec{r},\vec{\xi}) + k_0^2 G_0 (\vec{r},\vec{\xi}))\\
\nonumber &=& G_0 (\vec{r},\vec{\xi}) (-2 k_0^2 n_1(\vec{r}) + e^{- \Phi_0 (\vec{r},\vec{\xi})}\Phi_1 (\vec{r},\vec{\xi}) \delta (|\vec{r}-\vec{\xi}|))\\ \nonumber &\ & - \Phi_1 (\vec{r},\vec{\xi}) \delta (|\vec{r}-\vec{\xi}|)\\
&=& -2 G_0 (\vec{r},\vec{\xi}) k_0^2 n_1(\vec{r})
\end{eqnarray}

The solution to this linear differential equation with constant coefficients is obtained with the Green's function for free space propagation as
\begin{eqnarray}
\Phi_1 (\vec{r},\vec{\xi}) 
\nonumber &=& \frac{W_1 (\vec{r},\vec{\xi})}{G_0 (\vec{r},\vec{\xi})}\\
&=& \frac{1}{G_0 (\vec{r},\vec{\xi})} \int 2 G_0 (\vec{\zeta},\vec{\xi}) k_0^2 n_1(\vec{\zeta})G_0 (\vec{r},\vec{\zeta}) d\vec{\zeta}
\end{eqnarray}
where the integration is over the source volume denoted by primed coordinates, and $G_0$ is given by (\ref{eq:3.10}).

For simplification, we denote
\begin{equation} \label{eq:defK}
K(\vec{r}, \vec{\xi}, \vec{\zeta}) = \frac{G_0 (\vec{\zeta},\vec{\xi}) G_0 (\vec{r},\vec{\zeta})}{G_0 (\vec{r},\vec{\xi})}
\end{equation}
then the Green's function for (\ref{eq:1cm1}) to the second order is
\begin{eqnarray}  
G(\vec{r},\vec{\xi}) 
\nonumber &=& e^{\Phi(\vec{r},\vec{\xi})}\\
\label{eq:3.17} &=& G_0(\vec{r},\vec{\xi}) e^{\delta \int 2k_0^2 n_1(\vec{\zeta}) K(\vec{r}, \vec{\xi}, \vec{\zeta}) d \zeta}
\end{eqnarray}

\section{Reciprocity of the Green's Function}

Theorem\ \ \ \ \ The Green's function for the 3-D Helmholtz equation is symmetric. In other words, if the function $G$ Satisfies the Helmholtz equation
\begin{equation} \label{eq:3.18}
\Delta G(\vec{r},\vec{r}_0) + k_0^2 n^2(\vec{r}) G(\vec{r},\vec{r}_0) = -\delta(\vec{r}-\vec{r}_0)
\end{equation}
where $k_0$ is the wavenumber, $n(\vec{r})$ is the refractive index, and $\delta (\cdot)$ is the volume impulse function, then
\begin{equation}
G(\vec{r}_1,\vec{r}_0) = G(\vec{r}_0,\vec{r}_1)
\end{equation}

Proof: We shall assume that $G$ satisfies a homogeneous boundary condition at infinity, i.e.
\begin{equation} \label{eq:star}
|G(\vec{r},\vec{r}_0)| \to 0 \ \ \ \ \ as \ \ \ \ \  |\vec{r}-\vec{r}_0| \to \infty
\end{equation}

From equation (\ref{eq:3.18}) we have that
\begin{equation} \label{eq:3.19}
G(\vec{r},\vec{r}_1)\Delta G(\vec{r},\vec{r}_0) - G(\vec{r},\vec{r}_0) \Delta G(\vec{r},\vec{r}_1)= - [\delta(\vec{r}-\vec{r}_0) G(\vec{r},\vec{r}_1) - \delta(\vec{r}-\vec{r}_1) G(\vec{r},\vec{r}_0) ]
\end{equation}

Let $S$ denote the sphere of infinite radius centered at $r_0$, and $V$ denote the volume enclosed by $S$, then we may integrate (\ref{eq:3.19}) over $V$ and obtain
\begin{equation} \label{eq:3.20}
\int_V [G(\vec{r},\vec{r}_1)\Delta G(\vec{r},\vec{r}_1) - G(\vec{r},\vec{r}_0) \Delta G(\vec{r},\vec{r}_1)] dV = G(\vec{r}_1, \vec{r}_0) - G(\vec{r}_0, \vec{r}_1)
\end{equation}

Applying Green's theorem on the LHS of (\ref{eq:3.20}), we obtain
\begin{equation}
\int_S [G(\vec{r},\vec{r}_1) \nabla G(\vec{r},\vec{r}_0)- G(\vec{r},\vec{r}_0) \nabla G(\vec{r},\vec{r}_1)] \cdot dA = G(\vec{r}_1,\vec{r}_0) - G(\vec{r}_0,\vec{r}_1)
\end{equation}
and from the homogeneous boundary condition (\ref{eq:star}), we obtain
\begin{equation}
G(\vec{r}_1,\vec{r}_0) = G(\vec{r}_0,\vec{r}_1)
\end{equation}
End of proof.\\

It's easy to check that the Green's function (\ref{eq:3.17}) we have obtained by Rytov method satisfies symmetry property, since
\begin{equation}
G_0 (\vec{r},\vec{\xi})= \frac{e^{ik_0 |\vec{r}-\vec{\xi}|}}{4\pi |\vec{r}-\vec{\xi}|}
\end{equation}
is a symmetric function, thus $K(\vec{r},\vec{\xi},\vec{\zeta})=K(\vec{\xi},\vec{r},\vec{\zeta})$, where K is defined in (\ref{eq:defK})

\newpage
\pagestyle{myheadings} 
\markright{  \rm \normalsize CHAPTER 6. \hspace{0.5cm} 
  Parabolic Approximation}

\chapter{Parabolic Approximation}
\thispagestyle{myheadings}

\section{Solution of Helmholtz Equation by Parabolic Approximation}

As in the previous two chapters, the starting point of all propagation theory is the Helmholtz equation
\begin{equation} \label{eq:pp1}
\Delta u(\vec{r}) + k_0^2 n^2(\vec{r}) u(\vec{r})=0
\end{equation}
where $k_0$ is the wave number, and $n(\cdot)$ is the refractive index, for which we can assume the form
\begin{equation} \label{eq:ppa}
n(\vec{r}) = 1+\delta n_1(\vec{r})
\end{equation}

Next we choose a propagation direction: say, X-axis, and write
\begin{eqnarray}
\vec{r} = \vec{i} x + \vec{j} y + \vec{k} z\\
\rho = \vec{j} y + \vec{k} z
\end{eqnarray}
where $\vec{i}$, $\vec{j}$, $\vec{k}$ are orthonormal unit vectors in 3-D space, and $x$, $y$, $z$ are scalars. Then we let
\begin{equation}
u(x,\rho) = v(x, \rho) e^{i k x}
\end{equation}
and substituting in the Helmholtz equation (\ref{eq:pp1}) we obtain
\begin{equation} \label{eq:pp2}
\frac{\partial^2 v}{\partial x^2} +2 i k_0 \frac{\partial v}{\partial x} + k_0^2 (n^2 -1) v + \Delta_{\rho} v =0
\end{equation}
where $\Delta_{\rho} v = ( \frac{\partial^2}{\partial y^2} +\frac{\partial^2}{\partial z^2}) v$.\\
Equation (\ref{eq:pp2}) is equivalent to
\begin{equation}
\frac{\partial^2 v}{\partial x^2} +2 i k_0 \frac{\partial v}{\partial x} = \frac{\partial}{\partial x} (\frac{\partial v}{\partial x} + 2 i k_0 v)
\end{equation}

We assume that $v$ is slow-varying in $x$, so that
\begin{equation}
\frac{\partial v}{\partial x} + 2 i k_0 v \approx 2 i k_0 v
\end{equation}
or that the $\frac{\partial^2 v}{\partial x^2}$ term can be neglected. Then we have the parabolic approximation
\begin{equation} \label{eq:qq1}
2 i k_0 \frac{\partial v}{\partial x} + k_0^2 (n^2 -1) v + \Delta_{\rho} v =0
\end{equation}

In Rytov's method the solution is expressed as 
\begin{equation}
v(x, \rho)=exp[\Phi(x, \rho)]
\end{equation}
This leads to the following equation
\begin{equation} \label{eq:pp4}
2 i k_0 \frac{\partial \Phi}{\partial x} + \Delta_{\rho} \Phi + \nabla_{\rho} \Phi \cdot \nabla_{\rho} \Phi = - k_0^2 (n^2 -1)
\end{equation}
where $\nabla_{\rho} \Phi = (\frac{\partial \Phi}{\partial y}, \frac{\partial \Phi}{\partial z})$.\\
In the perturbation method, we assume that the solution can be written as a power series in $\delta$, i.e.
\begin{equation} \label{eq:pp3}
\Phi(x, \rho, \delta)=\Phi_0(x, \rho) + \delta \Phi_1(x, \rho)+ \delta^2 \Phi_2(x, \rho)+......
\end{equation}
By inserting (\ref{eq:ppa}) and (\ref{eq:pp3}) into (\ref{eq:pp4}), and by separating the equations based on the power of $\delta$, we obtain the equation for $\delta^0$ as
\begin{equation} \label{eq:pp5}
2 i k_0 \frac{\partial \Phi_0}{\partial x} + \Delta_{\rho} \Phi_0 + \nabla_{\rho} \Phi_0 \cdot \nabla_{\rho} \Phi_0 = 0
\end{equation}
Now we set 
\begin{equation}
v_0 (x,\rho) = e^{\Phi_0 (x, \rho)}
\end{equation}
then equation (\ref{eq:pp5}) is equivalent to

\begin{equation} \label{eq:pp6}
2 i k_0 \frac{\partial v_0}{\partial x} + \Delta_{\rho} v_0 = 0
\end{equation}

The equation for $\delta^1$ is 
\begin{equation} \label{eq:pp7}
2 i k_0 \frac{\partial \Phi_1}{\partial x} + \Delta_{\rho} \Phi_1 + 2 \nabla_{\rho} \Phi_1 \cdot \nabla_{\rho} \Phi_0 = -2 n_1 k_0^2
\end{equation}
If I set
\begin{equation} \label{eq:pp8}
W_1(x,\rho) = \Phi_1(x,\rho) v_0(x,\rho)
\end{equation}
Now by combining (\ref{eq:pp6}), (\ref{eq:pp7}) and (\ref{eq:pp8}), we have
\begin{eqnarray}
\nonumber &\ &2 i k_0 \frac{\partial W_1}{\partial x} + \Delta_{\rho} W_1\\
\nonumber &=&2 i k_0\frac{\partial}{\partial x}(\Phi_1 v_0) + \Delta_{\rho} (\Phi_1 v_0)\\
\nonumber &=&2 i k_0 \frac{\partial \Phi_1}{\partial x} v_0 + 2 i k_0 \Phi_1 \frac{\partial v_0}{\partial x} + v_0 \Delta_{\rho} \Phi_1 + 2 \nabla_{\rho}\Phi_1 \cdot \nabla_{\rho} v_0 + \Phi_1 \Delta_{\rho} v_0\\
&=&- 2 k_0^2 n_1 v_0
\end{eqnarray}

The solution to this equation is obtained with the Green's function for free space propagation as
\begin{eqnarray} \label{eq:pp10}
\nonumber \Phi_1(x, \rho) &=& \frac{W_1 (x, \rho)}{v_0 (x, \rho)}\\
&=& \frac{2k_0^2}{v_0 (x, \rho)} \int v_0(\sigma, \theta) n_1(\sigma, \theta) G_p (x, \rho; \sigma, \theta) d \sigma d \theta
\end{eqnarray}
where the integration is over the source volume, and $G_p$ is the Green's function for (\ref{eq:pp6}).\\
The Green's function has been obtained by Ishimaru [21] as
\begin{equation}
G_p (x, \rho; \sigma, \theta) = \frac{e^{\frac{ik|\rho - \theta|^2}{2(x-\sigma)}}}{4 \pi (x-\sigma)}
\end{equation}
Thus the perturbed field is
\begin{equation} \label{eq:pp12}
\Phi_1(x, \rho) = \frac{2k_0^2}{v_0 (x, \rho)} \int v_0(\sigma, \theta) n_1(\sigma, \theta) \frac{e^{\frac{ik|\rho - \theta|^2}{2(x-\sigma)}}}{4 \pi (x-\sigma)} d \sigma d \theta
\end{equation}

The solution to the 2nd order is then 
\begin{equation}
v(x,\rho) = v_0(x, \rho) e^{\delta \Phi_1(x, \rho)}
\end{equation}
and
\begin{equation}
u(x, \rho)=v(x, \rho) e^{ikx} = v_0(x, \rho) e^{\delta \Phi_1(x, \rho) + ikx}
\end{equation}

We specialize next to two cases.

\textbf{Case 1}. \textit{Plane Wave}\\
Here we take 
\begin{equation}
u_0(x, \rho)= e^{ikx}
\end{equation}
and hence
\begin{equation}
v_0(x, \rho) = 1
\end{equation}
since $u_0=v_0 e^{ikx}$. Then equation (\ref{eq:pp12}) yields
\begin{eqnarray}
\nonumber \Phi_1(x, \rho) &=& 2k_0^2 \int n_1(\sigma, \theta) \frac{e^{\frac{ik|\rho - \theta|^2}{2(x-\sigma)}}}{4 \pi (x-\sigma)} d \sigma d \theta\\
&=& \frac{k_0^2}{2 \pi} \int \frac{n_1(\sigma, \theta)}{x-\sigma} e^{\frac{ik|\rho - \theta|^2}{2(x-\sigma)}}d \sigma d \theta
\end{eqnarray}

\textbf{Case 2}. \textit{Beam Wave}\\
For the beam wave case
\begin{equation}
u_0(0, \rho) = exp( - \frac{\rho^2 \alpha k}{2})
\end{equation}
and hence by Ishimaru [21], we have
\begin{equation}
u_0(x, \rho)= \frac{e^{ikx}}{1+ix\alpha} exp(-\frac{\alpha k \rho^2}{2(1+ix\alpha)})
\end{equation} 
which is equivalent to
\begin{equation}
v_0(x, \rho)= \frac{1}{1+ix\alpha} exp(-\frac{\alpha k \rho^2}{2(1+ix\alpha)})
\end{equation} 
Hence
\begin{eqnarray}
\nonumber \Phi_1(x, \rho) &=& \frac{2k_0^2}{\frac{1}{1+ix\alpha} exp(-\frac{\alpha k \rho^2}{2(1+ix\alpha)})} \int \frac{1}{1+i\sigma\alpha} exp(-\frac{\alpha k \theta^2}{2(1+i\sigma\alpha)}) n_1(\sigma, \theta) \frac{e^{\frac{ik|\rho - \theta|^2}{2(x-\sigma)}}}{4 \pi (x-\sigma)} d \sigma d \theta\\
\nonumber &=&\frac{k_0^2}{2 \pi} (1+ix\alpha) exp(\frac{\alpha k \rho^2}{2(1+ix\alpha)})\\
\nonumber &\ & \int \frac{n_1(\sigma, \theta)}{(1+i\sigma\alpha)(x-\sigma)} exp(-\frac{\alpha k \theta^2}{2(1+i\sigma\alpha)}+\frac{ik|\rho - \theta|^2}{2(x-\sigma)})d \sigma d \theta\\
\end{eqnarray}

\section{Green's Function by Parabolic Approximation} 
The Green's function $G'(x, \rho; \sigma, \theta)$ for (\ref{eq:qq1}) satisfies the following equation
\begin{equation} \label{eq:qqq}
2 i k_0 \frac{\partial}{\partial x} G'(x, \rho; \sigma, \theta) + k_0^2 (n^2 -1) G'(x, \rho; \sigma, \theta) + \Delta_{\rho} G'(x, \rho; \sigma, \theta) = -\delta (x, \rho; \sigma, \theta)
\end{equation}
where $\delta(\cdot;\cdot)$ is the Dirac delta function.\\
In Rytov's method, the solution is expressed as
\begin{equation}
G'(x, \rho; \sigma, \theta)=exp(\Phi(x, \rho; \sigma, \theta))
\end{equation}
which leads to the following equation
\begin{equation} \label{eq:qq2}
2 i k_0 \frac{\partial \Phi}{\partial x} + \Delta_{\rho} \Phi + \nabla_{\rho} \Phi \cdot \nabla_{\rho} \Phi = - k_0^2 (n^2 -1) - e^{-\Phi} \delta (x, \rho; \sigma, \theta)
\end{equation}  
We assume that the solution to (\ref{eq:qq2}) can be written as a power series in $\delta$,
\begin{equation} \label{eq:qq3}
\Phi (x, \rho; \sigma, \theta; \delta)=\Phi_0 (x, \rho; \sigma, \theta) + \delta \Phi_1 (x, \rho; \sigma, \theta)+ \delta^2 \Phi_2 (x, \rho; \sigma, \theta) + \dots
\end{equation}
Notice that the refractive index assumes the form $n(\cdot)=1+\delta n_1(\cdot)$, the right hand side of equation (\ref{eq:qq2}) could be written as
\begin{equation}
RHS=-k_0^2 (2\delta n_1 +\delta^2 n_1^2)-e^{- \Phi_0 } (1-\delta \Phi_1 + O (\delta^2)) \delta (x, \rho; \sigma, \theta)
\end{equation}
where we have approximated $\Phi$ by the first two terms in (\ref{eq:qq3}), and applied Tayler's approximation on $exp(-\delta \Phi_1)$.
By inserting (\ref{eq:qq3}) into (\ref{eq:qq2}) and by separating the equations based on the power of $\delta$, we obtained the equation for $\delta^0$ as
\begin{equation} \label{eq:qq4}
2 i k_0 \frac{\partial \Phi_0}{\partial x} + \Delta_{\rho} \Phi_0 + \nabla_{\rho} \Phi_0 \cdot \nabla_{\rho} \Phi_0 =  - e^{-\Phi_0} \delta (x, \rho; \sigma, \theta)
\end{equation}
Now we set
\begin{equation}
G_0(x, \rho; \sigma, \theta)=exp(\Phi_0 (x, \rho; \sigma, \theta))
\end{equation}
Then equation (\ref{eq:qq4}) is equivalent to 
\begin{equation} \label{eq:qq5}
2 i k_0 \frac{\partial}{\partial x} G_0(x, \rho; \sigma, \theta) + \Delta_{\rho} G_0(x, \rho; \sigma, \theta) = -\delta (x, \rho; \sigma, \theta)
\end{equation}
As we know, the Green's function for free-space propagation has been obtained as
\begin{equation} \label{eq:qq7}
G_0 (x, \rho; \sigma, \theta) = \frac{e^{\frac{ik|\rho - \theta|^2}{2(x-\sigma)}}}{4 \pi (x-\sigma)}
\end{equation}
The equation for $\delta^1$ is
\begin{equation} \label{eq:qq6}
2 i k_0 \frac{\partial \Phi_1}{\partial x} + \Delta_{\rho} \Phi_1 + 2 \nabla_{\rho} \Phi_1 \cdot \nabla_{\rho} \Phi_0 = - 2 k_0^2 n_1 + \Phi_1 e^{-\Phi_0} \delta (x, \rho; \sigma, \theta)
\end{equation}
If I set
\begin{equation}
W_1(x, \rho; \sigma, \theta)= \Phi_1(x, \rho; \sigma, \theta) G_0(x, \rho; \sigma, \theta)
\end{equation}
Then we obtain
\begin{eqnarray}
\nonumber &\ &2 i k_0 \frac{\partial W_1}{\partial x} + \Delta_{\rho} W_1\\
\nonumber &=&2 i k_0\frac{\partial}{\partial x}(\Phi_1 G_0) + \Delta_{\rho} (\Phi_1 G_0)\\
\nonumber &=&2 i k_0 \frac{\partial \Phi_1}{\partial x} G_0 + 2 i k_0 \Phi_1 \frac{\partial G_0}{\partial x} + G_0 \Delta_{\rho} \Phi_1 + 2 \nabla_{\rho}\Phi_1 \cdot \nabla_{\rho} G_0 + \Phi_1 \Delta_{\rho} G_0\\
&=&- 2 k_0^2 n_1 G_0
\end{eqnarray}
where we have applied (\ref{eq:qq5}) and (\ref{eq:qq6}).
The solution to this linear differential equation is obtained with the Green's function for free space propagation as
\begin{eqnarray} \label{eq:pp10}
\nonumber \Phi_1(x, \rho; \sigma, \theta) &=& \frac{W_1 (x, \rho; \sigma, \theta)}{G_0 (x, \rho; \sigma, \theta)}\\
\nonumber &=& \frac{2k_0^2}{G_0 (x, \rho; \sigma, \theta)} \int G_0(\tau, \gamma; \sigma, \theta) n_1(\tau, \gamma) G_0 (x, \rho; \tau, \gamma) d \tau d \gamma\\
\end{eqnarray}
where the integration is over the source volume, and $G_0$ is given by (\ref{eq:qq7}).\\
For simplification, we denote
\begin{equation}
K'(x,\rho;\sigma, \theta; \tau, \gamma) = \frac{G_0(\tau, \gamma; \sigma, \theta)G_0 (x, \rho; \tau, \gamma)}{G_0 (x, \rho; \sigma, \theta)}
\end{equation}
then the Green's function for (\ref{eq:qqq}) is
\begin{eqnarray}
\nonumber G'(x, \rho; \sigma, \theta) &=& exp(\Phi(x, \rho; \sigma, \theta))\\
\nonumber &=& G_0(x, \rho; \sigma, \theta) exp(\delta \int 2 k_0^2 n_1(\tau, \gamma) K'(x,\rho;\sigma, \theta; \tau, \gamma) d \tau d \gamma) \label{eq:qq10}\\
\end{eqnarray}
It's easy to see that the Green's function (\ref{eq:qq10}) we have obtained by using parabolic approximation satisfies the symmetry property.

\newpage
\pagestyle{myheadings} 
\markright{  \rm \normalsize CHAPTER 7. \hspace{0.5cm} 
  Time Reversal Application}

\chapter{Time Reversal Application}
\thispagestyle{myheadings}

An important phenomenon in wave propagation in non-homogeneous media is related to time reversal of the wave field. In the time reversal procedure, the wave received by an antenna (receiver-emitter) array is recorded and then re-emitted into the medium time reversed, that is , the tails of the recorded signals are sent first. The time-reversal procedure is equivalent to phase conjugation on the spatial component of the time-harmonic wave field.\\

\begin{figure}[h]
\begin{center}
\includegraphics[width=12cm]{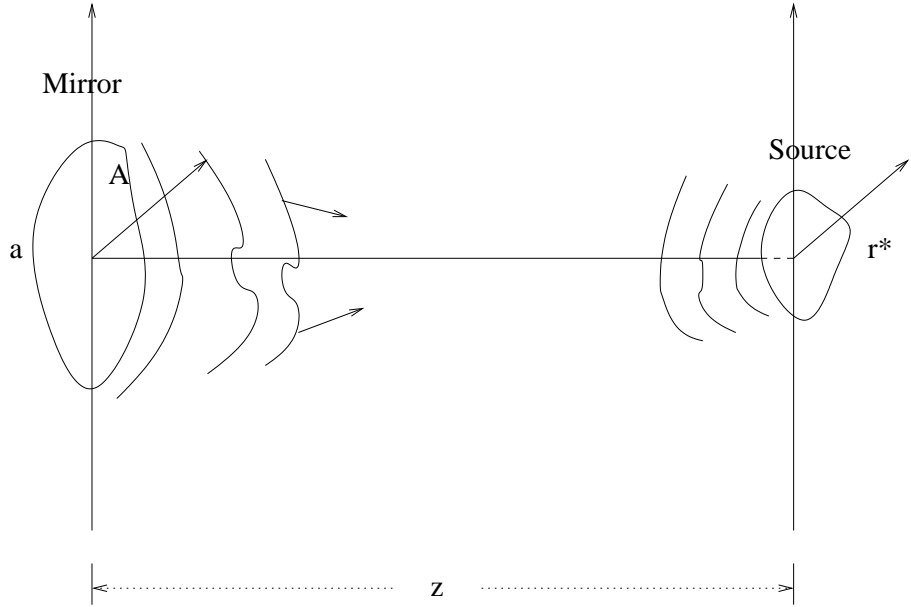}
\caption{The time reversal procedure. A pulse is emitted from a source with central wave length $\lambda_0$. The transmitted field is recorded, stored and time reversed at the mirror of size $a$ at $z$ distance away, and then sent back toward the source. It refocuses on the spot of size $r^\ast$. } \label{fig:5.1}
\end{center}
\end{figure}

A cartoon of a time reversal experiment is given in Figure (\ref{fig:5.1})

\newpage

To simplify the problem, let the phase-conjugated mirror be located at $z=0$, and the source at the parallel plane $z$-distance away. The aperture function of the mirror is assumed to be the indicator function $\chi_A$, where the set $A$ represents the physical boundary of the mirror.\\

Let $G_H (\vec{r}_0,\vec{r}_1)$ be the Green's function for the 3-D Helmholtz equation in random media, then $G (\vec{r}_0,\vec{r}_1)$ defined in (\ref{eq:3.17}) is an approximation.\\

In section 3.2, we stated and proved that $G (\vec{r}_0,\vec{r}_1)$ satisfies the symmetry property
\begin{equation}
G (\vec{r}_0,\vec{r}_1) = G (\vec{r}_1,\vec{r}_0)
\end{equation}

The wave field $\Psi_m$ received at the mirror is given by the super-position integral
\begin{equation}
\Psi_m (\vec{r}_m) = \chi_A (\vec{r}_m) \int G (\vec{r}_m,\vec{r}_s) \Psi_0 (\vec{r}_s) d \vec{r}_s
\end{equation}
where the integration is over the source volume denoted by primed coordinates, and $\Psi_0 (\vec{r}_s)$ is the input field.\\

After phase conjugation and back-propagation we have at the source plane the wave field.
\begin{equation} \label{eq:4.2}
\Psi^B (\vec{r}) = \int G (\vec{r},\vec{r}_m) \overline{G (\vec{r}_m,\vec{r}_s)} \chi_A (\vec{r}_m) \overline{\Psi_0 (\vec{r}_s)} d \vec{r}_m d \vec{r}_s
\end{equation}

By inserting (\ref{eq:3.17}) into (\ref{eq:4.2}), we obtain
\begin{eqnarray}
\Psi^B (\vec{r}) \nonumber &=& \int G_0 (\vec{r},\vec{r}_m) \overline{G_0 (\vec{r}_m,\vec{r}_s)} e^{\delta \int 2 k_0^2 n_1 (\vec{\zeta}) K (\vec{r}, \vec{r}_m, \vec{\zeta}) d \vec{\zeta}} \\
\nonumber &\ & \cdot e^{\delta \int 2 k_0^2 n_1 (\vec{\zeta}) \overline{K (\vec{r}_m, \vec{r}_s, \vec{\zeta})} d \vec{\zeta}} \chi_A 
(\vec{r}_m) \overline{\Psi_0 (\vec{r}_s)} d \vec{r}_m d \vec{r}_s\\
\nonumber &=&\int G_0 (\vec{r},\vec{r}_m) \overline{G_0 (\vec{r}_m,\vec{r}_s)} e^{\delta \int 2 k_0^2 n_1 (\vec{\zeta}) (K (\vec{r},\vec{r}_m, \vec{\zeta}) + \overline{K (\vec{r}_m, \vec{r}_s, \vec{\zeta})}) d \vec{\zeta}} \\
\label{eq:psi} &\ & \cdot \chi_A (\vec{r}_m) \overline{\Psi_0 (\vec{r}_s)} d \vec{r}_m d \vec{r}_s
\end{eqnarray}
where $G_0 (\vec{r},\vec{\xi})= \frac{e^{ik_0 |\vec{r}-\vec{\xi}|}}{4\pi |\vec{r}-\vec{\xi}|}$, and $ K(\vec{r}, \vec{\xi}, \vec{\zeta}) = \frac{G_0 (\vec{\zeta},\vec{\xi}) G_0 (\vec{r},\vec{\zeta})}{G_0 (\vec{r},\vec{\xi})}$.

Let $E(\cdot)$ denote the ensemble average, we can evaluate the field received at the source plane by
\begin{eqnarray}
E (\Psi^B (\vec{r}, n_1))
\nonumber &=& E( \int G_0 (\vec{r},\vec{r}_m) \overline{G_0 (\vec{r}_m,\vec{r}_s)} e^{\delta \int 2 k_0^2 n_1 (\vec{\zeta}) (K (\vec{r},\vec{r}_m, \vec{\zeta}) + \overline{K (\vec{r}_m, \vec{r}_s, \vec{\zeta})}) d \vec{\zeta}} \\
&\ & \cdot \chi_A (\vec{r}_m) \overline{\Psi_0 (\vec{r}_s)} d \vec{r}_m d \vec{r}_s)\\
\nonumber &=& \int G_0 (\vec{r},\vec{r}_m) \overline{G_0 (\vec{r}_m,\vec{r}_s)} E (e^{\delta \int 2 k_0^2 n_1 (\vec{\zeta}) (K (\vec{r},\vec{r}_m, \vec{\zeta}) + \overline{K (\vec{r}_m, \vec{r}_s, \vec{\zeta})}) d \vec{\zeta}}) \\
&\ & \cdot \chi_A (\vec{r}_m) \overline{\Psi_0 (\vec{r}_s)} d \vec{r}_m d \vec{r}_s\\
\nonumber &=& \int G_0 (\vec{r},\vec{r}_m) \overline{G_0 (\vec{r}_m,\vec{r}_s)}\\
\nonumber &\ & \cdot e^{4 \delta^2 k_0^4 \int E (n_1 (\vec{\zeta}_1) n_1 (\vec{\zeta}_2)) (K (\vec{r},\vec{r}_m, \vec{\zeta}_1) + \overline{K (\vec{r}_m, \vec{r}_s, \vec{\zeta}_1)}) (K (\vec{r},\vec{r}_m, \vec{\zeta}_2) + \overline{K (\vec{r}_m, \vec{r}_s, \vec{\zeta}_2)}) d \vec{\zeta}_1 d \vec{\zeta}_2} \\
&\ & \cdot \chi_A (\vec{r}_m) \overline{\Psi_0 (\vec{r}_s)} d \vec{r}_m d \vec{r}_s
\end{eqnarray}
where $E (n_1 (\vec{\zeta}_1) n_1 (\vec{\zeta}_2))$ is measured by experiments.\\

Experimentally, it is found that when the wave amplitude $\Psi^B$ is recorded, either photographically or by striking an observation screen, we record only the intensity $I=|\Psi^B|^2$. Hence, from (\ref{eq:psi}) we have
\begin{eqnarray}
I (\vec{r}) \nonumber &=& | \Psi^B(\vec{r}) |^2\\
\nonumber &=& \Psi^B(\vec{r}) \overline{\Psi^B(\vec{r})}\\
\nonumber &=& \int G_0 (\vec{r},\vec{r}_{m1}) \overline{G_0 (\vec{r}_{m1}, \vec{r}_{s1})} \chi_A (\vec{r}_{m1}) \overline{\Psi_0 (\vec{r}_{s1})} e^{\delta \int 2 k_0^2 n_1 (\vec{\zeta}) [K (\vec{r},\vec{r}_{m1}, \vec{\zeta}) + \overline{K (\vec{r}_{m1}, \vec{r}_{s1}, \vec{\zeta})}] d \vec{\zeta}}\\
\nonumber &\ & d \vec{r}_{m1} d \vec{r}_{s1} \int \overline{G_0 (\vec{r},\vec{r}_{m2})} G_0 (\vec{r}_{m2}, \vec{r}_{s2}) \chi_A (\vec{r}_{m2}) \Psi_0 (\vec{r}_{s2})\\
\nonumber &\ & \cdot e^{\delta \int 2 k_0^2 n_1 (\vec{\zeta}) [\overline{K (\vec{r},\vec{r}_{m2}, \vec{\zeta})} + K (\vec{r}_{m2}, \vec{r}_{s2}, \vec{\zeta})] d \vec{\zeta}} d \vec{r}_{m2} d \vec{r}_{s2}\\
\nonumber &=& \int G_0 (\vec{r},\vec{r}_{m1}) \overline{G_0 (\vec{r}_{m1}, \vec{r}_{s1})} \overline{G_0 (\vec{r},\vec{r}_{m2})} G_0 (\vec{r}_{m2}, \vec{r}_{s2}) \chi_A (\vec{r}_{m1}) \chi_A (\vec{r}_{m2}) \overline{\Psi_0 (\vec{r}_{s1})} \Psi_0 (\vec{r}_{s2})\\
\nonumber &\ & \cdot e^{\delta \int 2 k_0^2 n_1 (\vec{\zeta}) [K (\vec{r},\vec{r}_{m1}, \vec{\zeta}) + \overline{K (\vec{r}_{m1}, \vec{r}_{s1}, \vec{\zeta})} + \overline{K (\vec{r},\vec{r}_{m2}, \vec{\zeta})} + K (\vec{r}_{m2}, \vec{r}_{s2}, \vec{\zeta}) ] d \vec{\zeta}} d \vec{r}_{m1} d \vec{r}_{s1} d \vec{r}_{m2} d \vec{r}_{s2}\\
\end{eqnarray}
where $G_0 (\vec{r},\vec{\xi})= \frac{e^{ik_0 |\vec{r}-\vec{\xi}|}}{4\pi |\vec{r}-\vec{\xi}|}$, and $ K(\vec{r}, \vec{\xi}, \vec{\zeta}) = \frac{G_0 (\vec{\zeta},\vec{\xi}) G_0 (\vec{r},\vec{\zeta})}{G_0 (\vec{r},\vec{\xi})}$.\\

By evaluating the ensemble average of the intensity, we obtain
\begin{eqnarray}
E ( I(\vec{r}, n_1)) \nonumber &=& E (\int G_0 (\vec{r},\vec{r}_{m1}) \overline{G_0 (\vec{r}_{m1}, \vec{r}_{s1})} \overline{G_0 (\vec{r},\vec{r}_{m2})} G_0 (\vec{r}_{m2}, \vec{r}_{s2}) \chi_A (\vec{r}_{m1}) \chi_A (\vec{r}_{m2})\\
\nonumber &\ & \cdot \overline{\Psi_0 (\vec{r}_{s1})} \Psi_0 (\vec{r}_{s2}) exp\{\delta \int 2 k_0^2 n_1 (\vec{\zeta}) [K (\vec{r},\vec{r}_{m1}, \vec{\zeta}) + \overline{K (\vec{r}_{m1}, \vec{r}_{s1}, \vec{\zeta})}\\
\nonumber &\ & + \overline{K (\vec{r},\vec{r}_{m2}, \vec{\zeta})} + K (\vec{r}_{m2}, \vec{r}_{s2}, \vec{\zeta}) ] d \vec{\zeta}\} d \vec{r}_{m1} d \vec{r}_{s1} d \vec{r}_{m2} d \vec{r}_{s2})\\
\nonumber &=& \int G_0 (\vec{r},\vec{r}_{m1}) \overline{G_0 (\vec{r}_{m1}, \vec{r}_{s1})} \overline{G_0 (\vec{r},\vec{r}_{m2})} G_0 (\vec{r}_{m2}, \vec{r}_{s2}) \chi_A (\vec{r}_{m1}) \chi_A (\vec{r}_{m2})\\
\nonumber &\ & \cdot \overline{\Psi_0 (\vec{r}_{s1})} \Psi_0 (\vec{r}_{s2}) E (exp\{\delta \int 2 k_0^2 n_1 (\vec{\zeta}) [K (\vec{r},\vec{r}_{m1}, \vec{\zeta}) + \overline{K (\vec{r}_{m1}, \vec{r}_{s1}, \vec{\zeta})}\\
\nonumber &\ & + \overline{K (\vec{r},\vec{r}_{m2}, \vec{\zeta})} + K (\vec{r}_{m2}, \vec{r}_{s2}, \vec{\zeta}) ] d \vec{\zeta}\}) d \vec{r}_{m1} d \vec{r}_{s1} d \vec{r}_{m2} d \vec{r}_{s2}\\
\nonumber &=& \int G_0 (\vec{r},\vec{r}_{m1}) \overline{G_0 (\vec{r}_{m1}, \vec{r}_{s1})} \overline{G_0 (\vec{r},\vec{r}_{m2})} G_0 (\vec{r}_{m2}, \vec{r}_{s2}) \chi_A (\vec{r}_{m1}) \chi_A (\vec{r}_{m2})\\
\nonumber &\ & \cdot \overline{\Psi_0 (\vec{r}_{s1})} \Psi_0 (\vec{r}_{s2}) exp\{4 \delta^2 k_0^4 \int E (n_1 (\vec{\zeta}_1) n_1 (\vec{\zeta}_2)) [K (\vec{r},\vec{r}_{m1}, \vec{\zeta}_1)\\
\nonumber &\ & + \overline{K (\vec{r}_{m1}, \vec{r}_{s1}, \vec{\zeta}_1)} + \overline{K (\vec{r},\vec{r}_{m2}, \vec{\zeta}_1)} + K (\vec{r}_{m2}, \vec{r}_{s2}, \vec{\zeta}_1) ] [K (\vec{r},\vec{r}_{m1}, \vec{\zeta}_2)\\
\nonumber &\ & + \overline{K (\vec{r}_{m1}, \vec{r}_{s1}, \vec{\zeta}_2)} + \overline{K (\vec{r},\vec{r}_{m2}, \vec{\zeta}_2)} + K (\vec{r}_{m2}, \vec{r}_{s2}, \vec{\zeta}_2) ] d \vec{\zeta}_1 d \vec{\zeta}_2\} \\
&\ & d \vec{r}_{m1} d \vec{r}_{s1} d \vec{r}_{m2} d \vec{r}_{s2}
\end{eqnarray}
where the data we need to evaluate $E (n_1 (\vec{\zeta}_1) n_1 (\vec{\zeta}_2))$ is measured by experiments.\\

To illuminate this method, let's consider one of the simplest cases. Assume that $\Psi_0(\cdot)$ and $\chi_A(\cdot)$ are Dirac delta functions, which means that both the source and the phase-conjugating mirror is a single concentration point. Thus we have
\begin{equation}
\Psi_m (\vec{r}_m) = \delta (\vec{r}_m) G (\vec{r}_m,\vec{r}_s)
\end{equation}
and
\begin{equation} \label{eq:oo1}
\Psi^B (\vec{r}_s) = G (\vec{r}_s,\vec{r}_m) \overline{G (\vec{r}_m,\vec{r}_s)}
\end{equation}
where $\vec{r}_s$ and $\vec{r}_m$ are position vectors for the source and the mirror respectively.\\
By inserting (\ref{eq:3.17}) into (\ref{eq:oo1}), we obtain
\begin{eqnarray}
\Psi^B (\vec{r}_s) \nonumber &=& G_0 (\vec{r}_s,\vec{r}_m) \overline{G_0 (\vec{r}_m,\vec{r}_s)} e^{\delta \int 2 k_0^2 n_1 (\vec{\zeta}) K (\vec{r}_s, \vec{r}_m, \vec{\zeta}) d \vec{\zeta}} e^{\delta \int 2 k_0^2 n_1 (\vec{\zeta}) \overline{K (\vec{r}_m, \vec{r}_s, \vec{\zeta})} d \vec{\zeta}}\\
&=& |G_0 (\vec{r}_s,\vec{r}_m)|^2 e^{\delta \int 4 k_0^2 n_1 (\vec{\zeta}) Re(K (\vec{r}_s,\vec{r}_m, \vec{\zeta})) d \vec{\zeta}} 
\end{eqnarray}
where $G_0 (\vec{r},\vec{\xi})= \frac{e^{ik_0 |\vec{r}-\vec{\xi}|}}{4\pi |\vec{r}-\vec{\xi}|}$, and $ K(\vec{r}, \vec{\xi}, \vec{\zeta}) = \frac{G_0 (\vec{\zeta},\vec{\xi}) G_0 (\vec{r},\vec{\zeta})}{G_0 (\vec{r},\vec{\xi})}$, and we have applied the symmetry property of $G_0$.\\
We can also obtain the intensity as
\begin{eqnarray}
I (\vec{r}_s) \nonumber &=& | \Psi^B(\vec{r}_s) |^2\\
\nonumber &=& \Psi^B(\vec{r}_s) \overline{\Psi^B(\vec{r}_s)}\\
&=& |G_0 (\vec{r}_s,\vec{r}_m)|^4 e^{\delta \int 8 k_0^2 n_1 (\vec{\zeta}) Re(K (\vec{r}_s,\vec{r}_m, \vec{\zeta})) d \vec{\zeta}} 
\end{eqnarray}
By evaluating the ensemble average, we obtain
\begin{eqnarray}
E ( I(\vec{r}_s, n_1)) \nonumber &=& E (|G_0 (\vec{r}_s,\vec{r}_m)|^4 e^{\delta \int 8 k_0^2 n_1 (\vec{\zeta}) Re(K (\vec{r}_s,\vec{r}_m, \vec{\zeta})) d \vec{\zeta}})\\
\nonumber &=& |G_0 (\vec{r}_s,\vec{r}_m)|^4 E(e^{\delta \int 8 k_0^2 n_1 (\vec{\zeta}) Re(K (\vec{r}_s,\vec{r}_m, \vec{\zeta})) d \vec{\zeta}})\\
\nonumber &=& |G_0 (\vec{r}_s,\vec{r}_m)|^4 e^{64 \delta^2 k_0^4 \int E(n_1 (\vec{\zeta}_1)n_1 (\vec{\zeta}_2)) Re(K (\vec{r}_s,\vec{r}_m, \vec{\zeta}_1))Re(K (\vec{r}_s,\vec{r}_m, \vec{\zeta}_2))d \vec{\zeta}_1 d \vec{\zeta}_2}\\
&=& |G_0 (\vec{r}_s,\vec{r}_m)|^4 e^{64 \delta^2 \sigma^2 k_0^4 \int Re(K (\vec{r}_s,\vec{r}_m, \vec{\zeta}))^2 d \vec{\zeta}}
\end{eqnarray}
where we have taken the Gaussian expectation and assumed that $n_1(\cdot)$ is white noise, so that $E(n_1 (\vec{\zeta}_1)n_1 (\vec{\zeta}_2)) = \sigma^2 \delta(\vec{\zeta}_1-\vec{\zeta}_2)$.

\newpage
\pagestyle{myheadings} 
\markright{  \rm \normalsize CHAPTER 8. \hspace{0.5cm} 
  3-D Apodization Problem In Turbulence}

\chapter{3-D Apodization Problem In Turbulence}
\thispagestyle{myheadings}

In this chapter, we are concerned with the determination of the distribution of light over the exit pupil of an optical system required in order to achieve a desired distribution of illuminance over a given plane in the image field. This problem is known as the Apodization Problem, and we are particular interested in the determination of that amplitude distribution over a circular pupil which maximizes the fraction of the total energy that lies in a prescribed circle in the image plane.\\

\begin{figure}[h]
\begin{center}
\includegraphics[width=12cm]{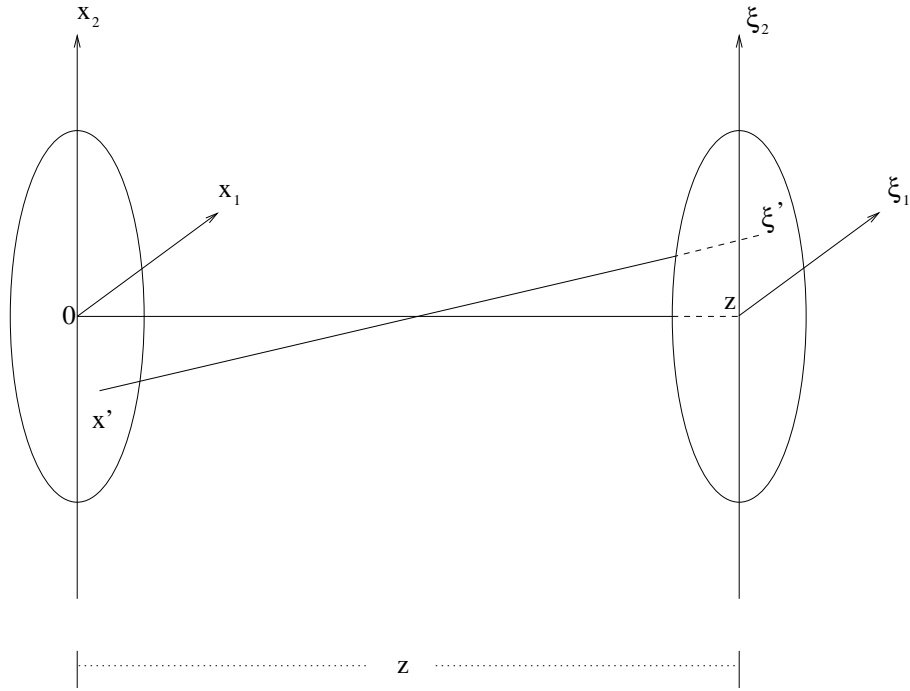}
\caption{A configuration of the Apodization Problem. } \label{fig:6.1}
\end{center}
\end{figure}

Let $\vec{x}=(x^1, x^2)$ be the radius vector in the plane of the exit pupil from the optical axis to an arbitrary point in that plane; let $\vec{\xi}=(\xi^1, \xi^2)$ be the radius vector in the image plane from the optical axis to a point in the image plane(Figure (\ref{fig:6.1})).
Then by the extended Huygens-Fresnal principle, the light amplitude $A(\vec{\xi})$ in the image plane is given by the superposition integral
\begin{equation} \label{eq:5.1}
A (\vec{\xi}) = \int_{|\vec{x}| \leq a} G_H (0, \vec{x}; z, \vec{\xi}) T(\vec{x}) d \vec{x}
\end{equation}
where $G_H (0, \vec{x}; z, \vec{\xi})$ is the atmospheric impulse response (Green's function), $T(\vec{x})$ is the light amplitude in the circular exit pupil of radius $a$, and $z$ is the distance from the pupil to the image plane.\\

Our apodization problem then requires finding the function $T(\vec{x})$ for which the ratio
\begin{equation} \label{eq:5.2}
\lambda = \frac{\int_{|\vec{\xi}| \leq b} |A (\vec{\xi})|^2 d \vec{\xi}} {\int_{R^2} |A (\vec{\xi})|^2 d \vec{\xi}}
\end{equation}
is a maximum.\\
Here $b$ is the radius of the circle in which the illuminance is to be maximally concentrated, and the circle is centered on the optical axis.\\

We assume the optical wave propagates through a thick slab of turbulence air, where Rytov's method could be applied to approximate the Helmholtz equation.\\
From chapter 4, we have derived an explicit formula for $G_H$ by Rytov's Approximation.
\begin{equation} \label{eq:5.3} 
G(\vec{r},\vec{\xi}) = G_0(\vec{r},\vec{\xi}) e^{\delta \int 2k_0^2 n_1(\vec{\zeta}) K(\vec{r}, \vec{\xi}, \vec{\zeta}) d \zeta}
\end{equation}
where
\begin{equation}
K(\vec{r}, \vec{\xi}, \vec{\zeta}) = \frac{G_0 (\vec{\zeta},\vec{\xi}) G_0 (\vec{r},\vec{\zeta})}{G_0 (\vec{r},\vec{\xi})}
\end{equation}
and
\begin{equation}
G_0 (\vec{r},\vec{\xi})= \frac{e^{ik_0 |\vec{r}-\vec{\xi}|}}{4\pi |\vec{r}-\vec{\xi}|}
\end{equation}

Denote the top of (\ref{eq:5.2}) by $M$, then $M$ is a function dependent of random variable $n_1(\vec{\zeta})$, i.e.
\begin{equation} \label{eq:5.4}
M(n_1) = \int_{|\vec{\xi}| \leq b} |A (\vec{\xi})|^2 d \vec{\xi}
\end{equation}

By inserting (\ref{eq:5.1}) and (\ref{eq:5.3}) into (\ref{eq:5.4}), we obtain
\begin{eqnarray}
M (n_1) \nonumber &=& \int_{|\vec{\xi}| \leq b} A (\vec{\xi}) \overline{A (\vec{\xi})} d \vec{\xi}\\
\nonumber &=& \int_{|\vec{\xi}| \leq b, |\vec{x}_2| \leq a, |\vec{x}_1| \leq a} G_H (0, \vec{x}_1; z, \vec{\xi}) T(\vec{x}_1) \overline{G_H (0, \vec{x}_2; z, \vec{\xi})} \overline{T(\vec{x}_2)} d \vec{x}_1 d \vec{x}_2 d \vec{\xi}\\
\nonumber &=& \int_{|\vec{\xi}| \leq b, |\vec{x}_2| \leq a, |\vec{x}_1| \leq a} T(\vec{x}_1) \overline{T(\vec{x}_2)} G_0 (0, \vec{x}_1; z, \vec{\xi}) \overline{G_0 (0, \vec{x}_2; z, \vec{\xi})}\\
\nonumber &\ & \cdot e^{\delta \int 2k_0^2 n_1(\vec{\zeta}) K[(0, \vec{x}_1), (z, \vec{\xi}), \vec{\zeta}] d \vec{\zeta}} e^{\delta \int 2k_0^2 n_1(\vec{\zeta}) \overline{K[(0, \vec{x}_2), (z, \vec{\xi}), \vec{\zeta}]} d \vec{\zeta}} d \vec{x}_1 d \vec{x}_2 d \vec{\xi}\\
\nonumber &=& \int_{|\vec{\xi}| \leq b, |\vec{x}_2| \leq a, |\vec{x}_1| \leq a} T(\vec{x}_1) \overline{T(\vec{x}_2)} G_0 (0, \vec{x}_1; z, \vec{\xi}) \overline{G_0 (0, \vec{x}_2; z, \vec{\xi})}\\
&\ & \cdot e^{\delta \int 2k_0^2 n_1(\vec{\zeta}) \{K[(0, \vec{x}_1), (z, \vec{\xi}), \vec{\zeta}] + \overline{K[(0, \vec{x}_2), (z, \vec{\xi}), \vec{\zeta}]} \} d \vec{\zeta}} d \vec{x}_1 d \vec{x}_2 d \vec{\xi}
\end{eqnarray}
where we have assumed $n_1(\vec{\zeta})$ to be real.\\

Let $E(\cdot)$ denote the ensemble average (expectation), we need to evaluate the statistics of the energy received at the source plane.
\begin{eqnarray}
E(M (n_1)) \nonumber &=& E (\int_{|\vec{\xi}| \leq b, |\vec{x}_2| \leq a, |\vec{x}_1| \leq a} T(\vec{x}_1) \overline{T(\vec{x}_2)} G_0 (0, \vec{x}_1; z, \vec{\xi}) \overline{G_0 (0, \vec{x}_2; z, \vec{\xi})}\\
\nonumber &\ & \cdot e^{\delta \int 2k_0^2 n_1(\vec{\zeta}) \{K[(0, \vec{x}_1), (z, \vec{\xi}), \vec{\zeta}] + \overline{K[(0, \vec{x}_2), (z, \vec{\xi}), \vec{\zeta}]} \} d \vec{\zeta}} d \vec{x}_1 d \vec{x}_2 d \vec{\xi} )\\
\nonumber &=& \int_{|\vec{\xi}| \leq b, |\vec{x}_2| \leq a, |\vec{x}_1| \leq a} T(\vec{x}_1) \overline{T(\vec{x}_2)} G_0 (0, \vec{x}_1; z, \vec{\xi}) \overline{G_0 (0, \vec{x}_2; z, \vec{\xi})}\\
\nonumber &\ & \cdot E (e^{\delta \int 2k_0^2 n_1(\vec{\zeta}) \{K[(0, \vec{x}_1), (z, \vec{\xi}), \vec{\zeta}] + \overline{K[(0, \vec{x}_2), (z, \vec{\xi}), \vec{\zeta}]} \} d \vec{\zeta}} ) d \vec{x}_1 d \vec{x}_2 d \vec{\xi}\\
\nonumber &=& \int_{|\vec{\xi}| \leq b, |\vec{x}_2| \leq a, |\vec{x}_1| \leq a} T(\vec{x}_1) \overline{T(\vec{x}_2)} G_0 (0, \vec{x}_1; z, \vec{\xi}) \overline{G_0 (0, \vec{x}_2; z, \vec{\xi})} exp(4 \delta^2 k_0^4\\
\nonumber &\ &   \int E( n_1(\vec{\zeta}_1) n_1(\vec{\zeta}_2) ) \{K[(0, \vec{x}_1), (z, \vec{\xi}), \vec{\zeta}_1] + \overline{K[(0, \vec{x}_2), (z, \vec{\xi}), \vec{\zeta}_1]} \}\\ 
\label{eq:aaa} &\ & \cdot \{ K[(0, \vec{x}_1), (z, \vec{\xi}), \vec{\zeta}_2] + \overline{K[(0, \vec{x}_2), (z, \vec{\xi}), \vec{\zeta}_2]} \}d \vec{\zeta}_1 d \vec{\zeta}_2) d \vec{x}_1 d \vec{x}_2 d \vec{\xi}
\end{eqnarray}

Assuming $n_1(\cdot)$ is white noise, i.e. $E(n_1 (\vec{\zeta}_1)n_1 (\vec{\zeta}_2)) = \sigma^2 \delta(\vec{\zeta}_1-\vec{\zeta}_2)$, we can simplify (\ref{eq:aaa}) as
\begin{eqnarray}
E(M (n_1)) \nonumber &=& \int_{|\vec{\xi}| \leq b, |\vec{x}_2| \leq a, |\vec{x}_1| \leq a} T(\vec{x}_1) \overline{T(\vec{x}_2)} G_0 (0, \vec{x}_1; z, \vec{\xi}) \overline{G_0 (0, \vec{x}_2; z, \vec{\xi})}\\
\nonumber &\ & \cdot e^{4 \delta^2 \sigma^2 k_0^4 \int (K[(0, \vec{x}_1), (z, \vec{\xi}), \vec{\zeta}] + \overline{K[(0, \vec{x}_2), (z, \vec{\xi}), \vec{\zeta}]} )^2 d \vec{\zeta}} d \vec{x}_1 d \vec{x}_2 d \vec{\xi}\\
\nonumber &=& \int_{|\vec{\xi}| \leq b, |\vec{x}_2| \leq a, |\vec{x}_1| \leq a} T(\vec{x}_1) \overline{T(\vec{x}_2)} G_0 (0, \vec{x}_1; z, \vec{\xi}) \overline{G_0 (0, \vec{x}_2; z, \vec{\xi})}\\
\nonumber &\ & \cdot e^{4 \delta^2 \sigma^2 k_0^4 \int (K[(0, \vec{x}_1), (z, \vec{\xi}), \vec{\zeta}]^2 + 2 K[(0, \vec{x}_1), (z, \vec{\xi}), \vec{\zeta}] \overline{K[(0, \vec{x}_2), (z, \vec{\xi}), \vec{\zeta}]} +\overline{K[(0, \vec{x}_2), (z, \vec{\xi}), \vec{\zeta}]}^2 ) d \vec{\zeta}}\\
\nonumber &\ & \cdot d \vec{x}_1 d \vec{x}_2 d \vec{\xi}\\
\nonumber &=& \int_{|\vec{\xi}| \leq b, |\vec{x}_2| \leq a, |\vec{x}_1| \leq a} T(\vec{x}_1) G_0 (0, \vec{x}_1; z, \vec{\xi}) e^{4 \delta^2 \sigma^2 k_0^4 \int K[(0, \vec{x}_1), (z, \vec{\xi}), \vec{\zeta}]^2 d \vec{\zeta}}\\
\nonumber &\ & \cdot \overline{T(\vec{x}_2)} \overline{G_0 (0, \vec{x}_2; z, \vec{\xi})} e^{4 \delta^2 \sigma^2 k_0^4 \int \overline{K[(0, \vec{x}_2), (z, \vec{\xi}), \vec{\zeta}]}^2 d \vec{\zeta}} \\
\label{eq:wn} &\ & \cdot e^{8 \delta^2 \sigma^2 k_0^4 \int K[(0, \vec{x}_1), (z, \vec{\xi}), \vec{\zeta}] \overline{K[(0, \vec{x}_2), (z, \vec{\xi}), \vec{\zeta}]} d \vec{\zeta}} d \vec{x}_1 d \vec{x}_2 d \vec{\xi}
\end{eqnarray}

Since $\delta<<1$, we can rewrite the last term inside the integral of (\ref{eq:wn}) by Tayler's expansion as
\begin{eqnarray}
\nonumber &\ &e^{8 \delta^2 \sigma^2 k_0^4 \int K[(0, \vec{x}_1), (z, \vec{\xi}), \vec{\zeta}] \overline{K[(0, \vec{x}_2), (z, \vec{\xi}), \vec{\zeta}]} d \vec{\zeta}} \\
&=& 1 + 8 \delta^2 \sigma^2 k_0^4 \int K[(0, \vec{x}_1), (z, \vec{\xi}), \vec{\zeta}] \overline{K[(0, \vec{x}_2), (z, \vec{\xi}), \vec{\zeta}]} d \vec{\zeta} + O(\delta^4)
\end{eqnarray}

Denote 
\begin{equation} \label{eq:bbc}
D(\vec{x}_1, \vec{\xi}) := T(\vec{x}_1) G_0 (0, \vec{x}_1; z, \vec{\xi}) e^{4 \delta^2 \sigma^2 k_0^4 \int K[(0, \vec{x}_1), (z, \vec{\xi}), \vec{\zeta}]^2 d \vec{\zeta}}
\end{equation}
Then (\ref{eq:wn}) is approximate by
\begin{eqnarray}
E(M (n_1)) \nonumber &\approx& \int_{|\vec{\xi}| \leq b, |\vec{x}_2| \leq a, |\vec{x}_1| \leq a} D(\vec{x}_1, \vec{\xi}) \overline{D(\vec{x}_2, \vec{\xi})} (1 + 8 \delta^2 \sigma^2 k_0^4 \int K[(0, \vec{x}_1), (z, \vec{\xi}), \vec{\zeta}]\\
\nonumber &\ & \cdot \overline{K[(0, \vec{x}_2), (z, \vec{\xi}), \vec{\zeta}]} d \vec{\zeta}) d \vec{x}_1 d \vec{x}_2 d \vec{\xi}\\
\nonumber &=& \int_{|\vec{\xi}| \leq b} |\int_{|\vec{x}| \leq a} D(\vec{x}, \vec{\xi}) d \vec{x}|^2 d \vec{\xi}\\
&+& 8 \delta^2 \sigma^2 k_0^4 \int_{|\vec{\xi}| \leq b} |\int_{|\vec{x}| \leq a} D(\vec{x}, \vec{\xi}) K[(0, \vec{x}), (z, \vec{\xi}), \vec{\zeta}] d \vec{x}|^2 d \vec{\zeta} d \vec{\xi}
\end{eqnarray}

If the air-density inhomogeneity is extremely small, which means $\delta^2 k_0^4<<1$, then the second term is negligible compared to the first one. Consequently we have
\begin{equation} \label{eq:ccb}
E(M (n_1)) =  \int_{|\vec{\xi}| \leq b} |\int_{|\vec{x}| \leq a} D(\vec{x}, \vec{\xi}) d \vec{x}|^2 d \vec{\xi}
\end{equation}

This result can be obtained equivalently by taking the expectation of the light amplitude $A (\vec{\xi})$ first, and then by evaluating the energy concentrated in the image circle from the expected light amplitude $E(A (\vec{\xi}))$, since
\begin{eqnarray}
\nonumber &\ & E(A (\vec{\xi}))\\
\nonumber &=& E (\int_{|\vec{x}| \leq a} G_H (0, \vec{x}; z, \vec{\xi}) T(\vec{x}) d \vec{x})\\
\nonumber &=& E(\int_{|\vec{x}| \leq a} G_0 (0, \vec{x}; z, \vec{\xi}) e^{\delta \int 2k_0^2 n_1(\vec{\zeta}) K[(0, \vec{x}), (z, \vec{\xi}), \vec{\zeta}] d \vec{\zeta}} T(\vec{x}) d \vec{x})\\
\nonumber &=& \int_{|\vec{x}| \leq a} T(\vec{x}) G_0 (0, \vec{x}; z, \vec{\xi}) E(e^{\delta \int 2k_0^2 n_1(\vec{\zeta}) K[(0, \vec{x}), (z, \vec{\xi}), \vec{\zeta}] d \vec{\zeta}}) d \vec{x}\\
\nonumber &=& \int_{|\vec{x}| \leq a} T(\vec{x}) G_0 (0, \vec{x}; z, \vec{\xi}) e^{4 \delta^2 k_0^4 \int E( n_1(\vec{\zeta}_1) n_1(\vec{\zeta}_2) ) K[(0, \vec{x}), (z, \vec{\xi}), \vec{\zeta}_1] K[(0, \vec{x}), (z, \vec{\xi}), \vec{\zeta}_2] d \vec{\zeta}_1 d \vec{\zeta}_2} d \vec{x}\\
\nonumber &=& \int_{|\vec{x}| \leq a} T(\vec{x}) G_0 (0, \vec{x}; z, \vec{\xi}) e^{4 \delta^2 \sigma^2 k_0^4 \int K[(0, \vec{x}), (z, \vec{\xi}), \vec{\zeta}]^2 d \vec{\zeta}} d \vec{x}\\
&=& \int_{|\vec{x}| \leq a} D(\vec{x}, \vec{\xi}) d \vec{x}
\end{eqnarray}
where
\begin{equation}
D(\vec{x}, \vec{\xi}) = T(\vec{x}) G_0 (0, \vec{x}; z, \vec{\xi}) e^{4 \delta^2 \sigma^2 k_0^4 \int K[(0, \vec{x}), (z, \vec{\xi}), \vec{\zeta}]^2 d \vec{\zeta}}
\end{equation}
which is the same as (\ref{eq:bbc}).\\
Therefore the energy received in the image circle is
\begin{eqnarray}
M' \nonumber &=& \int_{|\vec{\xi}| \leq b} |E(A (\vec{\xi}))|^2 d \vec{\xi}\\
&=& \int_{|\vec{\xi}| \leq b} |\int_{|\vec{x}| \leq a} D(\vec{x}, \vec{\xi}) d \vec{x}|^2 d \vec{\xi}
\end{eqnarray}
which is identical to (\ref{eq:ccb}).\\

Furthermore we can assume that the energy is conserved through wave propagation. Therefore by similar calculation, we obtain the following equation for the the bottom of (\ref{eq:5.2})
\begin{equation}
\int_{R^2} |A (\vec{\xi})|^2 d \vec{\xi} = \int_{R^2} |\int_{|\vec{x}| \leq a} D(\vec{x}, \vec{\xi}) d \vec{x}|^2 d \vec{\xi}
\end{equation}
Thus our original apodization problem is reduced to the one that requires finding the function $T(\vec{x})$ for which the ratio
\begin{equation} \label{eq:ene}
\lambda ' = \frac{\int_{|\vec{\xi}| \leq b} |\int_{|\vec{x}| \leq a} D(\vec{x}, \vec{\xi}) d \vec{x}|^2 d \vec{\xi}}
{\int_{R^2} |\int_{|\vec{x}| \leq a} D(\vec{x}, \vec{\xi}) d \vec{x}|^2 d \vec{\xi}}
\end{equation}
is the maximum, where $D(\vec{x}, \vec{\xi}) := T(\vec{x}) G_0 (0, \vec{x}; z, \vec{\xi}) e^{4 \delta^2 \sigma^2 k_0^4 \int K[(0, \vec{x}), (z, \vec{\xi}), \vec{\zeta}]^2 d \vec{\zeta}}$, and $K(\vec{r}, \vec{\xi}, \vec{\zeta}) = \frac{G_0 (\vec{\zeta},\vec{\xi}) G_0 (\vec{r},\vec{\zeta})}{G_0 (\vec{r},\vec{\xi})}$.\\

If A is of total energy $E$, since the energy is concerved, we have
\begin{equation}
E= \int_{R^2} |\int_{|\vec{x}| \leq a} D(\vec{x}, \vec{\xi}) d \vec{x}|^2 d \vec{\xi} = \int_{|\vec{x}| \leq a} |T(\vec{x})|^2 d \vec{x}
\end{equation}
whereas the energy in the image circle is
\begin{eqnarray}
\nonumber &\ & \int_{|\vec{\xi}| \leq b} |\int_{|\vec{x}| \leq a} D(\vec{x}, \vec{\xi}) d \vec{x}|^2 d \vec{\xi}\\
\nonumber &=& \int_{|\vec{\xi}| \leq b} \int_{|\vec{x}| \leq a} \int_{|\vec{y}| \leq a} D(\vec{x}, \vec{\xi}) \overline{D(\vec{y}, \vec{\xi})} d \vec{x} d \vec{y} d \vec{\xi}\\
&=& \int_{|\vec{x}| \leq a} \int_{|\vec{y}| \leq a} K_s(\vec{x}, \vec{y}) T(\vec{x}) \overline{T(\vec{y})} d \vec{x} d \vec{y}
\end{eqnarray}
where
\begin{equation}
K_s (\vec{x}, \vec{y}) = \int_{|\vec{\xi}| \leq b} G_0 (0, \vec{x}; z, \vec{\xi}) \overline{G_0 (0, \vec{y}; z, \vec{\xi})} e^{4 \delta^2  \sigma^2 k_0^4 \int (K[(0, \vec{x}), (z, \vec{\xi}), \vec{\zeta}]^2 + \overline{K[(0, \vec{y}), (z, \vec{\xi}), \vec{\zeta}]}^2) d \vec{\zeta}} d \vec{\xi}
\end{equation}

Our apodization problem therefore requires finding the function $T(\vec{x})$ for which the ratio
\begin{equation} \label{eq:ene}
\lambda ' = \frac{\int_{|\vec{x}| \leq a} \int_{|\vec{y}| \leq a} K_s(\vec{x}, \vec{y}) T(\vec{x}) \overline{T(\vec{y})} d \vec{x} d \vec{y}}
{\int_{|\vec{x}| \leq a} |T(\vec{x})|^2 d \vec{x}}
\end{equation}
is the maximum. This maximum is $\lambda_0$, the largest eigenvalue of the integral equation
\begin{equation} \label{eq:ii1}
\lambda ' \Psi (\vec{x}) = \int_{|\vec{y}| \leq a} K_s (\vec{x}, \vec{y}) \Psi (\vec{y}) d \vec{y}, \ \ \ \ |\vec{x}| \leq a
\end{equation}
with kernel
\begin{equation} \label{eq:ii4}
K_s (\vec{x}, \vec{y}) = \int_{|\vec{\xi}| \leq b} G_0 (0, \vec{x}; z, \vec{\xi}) \overline{G_0 (0, \vec{y}; z, \vec{\xi})} e^{4 \delta^2 \sigma^2 k_0^4 \int (K[(0, \vec{x}), (z, \vec{\xi}), \vec{\zeta}]^2 + \overline{K[(0, \vec{y}), (z, \vec{\xi}), \vec{\zeta}]}^2) d \vec{\zeta}} d \vec{\xi}
\end{equation}
Somewhat simpler than (\ref{eq:ii1}) is the integral equation
\begin{equation} \label{eq:ii2}
\alpha \Psi (\vec{x}) = \int_{|\vec{\eta}| \leq a} G_0 (0, \vec{x}; z, \vec{\eta}) e^{4 \delta^2 \sigma^2 k_0^4 \int K[(0, \vec{x}), (z, \vec{\eta}), \vec{\zeta}]^2 d \vec{\zeta}} \Psi (\vec{\eta}) d \vec{\eta}
\end{equation}
Whe shall show below that the solution of this equation is completely equivalent to the solution of (\ref{eq:ii1}).\\

From the symmetry of the input domain, it follows that if $\Psi (\vec{x})$ is a solution of (\ref{eq:ii2}), so also is $\Psi (-\vec{x})$, so that both $\Psi_e (\vec{x}) = \Psi (\vec{x}) + \Psi (-\vec{x})$ and $\Psi_o (\vec{x}) = \Psi (\vec{x}) - \Psi (-\vec{x})$ are solutions as well. The eigenfunctions of (\ref{eq:ii2}) can be chosen to be either even or odd functions of $\vec{x}$.\\
The complex conjugate of (\ref{eq:ii2}) is
\begin{equation} \label{eq:jjk}
\overline{\alpha} \overline{\Psi (\vec{x})} = \int_{|\vec{\eta}| \leq a} \overline{G_0 (0, \vec{x}; z, \vec{\eta})} e^{4 \delta^2 \sigma^2 k_0^4 \int \overline{K[(0, \vec{x}), (z, \vec{\eta}), \vec{\zeta}]}^2 d \vec{\zeta}} \overline{\Psi (\vec{\eta})} d \vec{\eta}
\end{equation}
Multiply (\ref{eq:ii2}) by $\overline{\Psi (\vec{x})}$ and integrate over $|\vec{x}| \leq a$. Multiply (\ref{eq:jjk}) by $\Psi (\vec{x})$ and integrate over $|\vec{x}| \leq a$. Combining these equations, we find on using the symmetry of the integral domain that
\begin{eqnarray}
\nonumber &\ & (\alpha \pm \overline{\alpha}) \int_{|\vec{x}| \leq a} \Psi (\vec{x}) \overline{\Psi (\vec{x})} d \vec{x}\\
\nonumber &=& \int_{|\vec{\eta}| \leq a, |\vec{x}| \leq a} G_0 (0, \vec{x}; z, \vec{\eta}) e^{4 \delta^2 \sigma^2 k_0^4 \int K[(0, \vec{x}), (z, \vec{\eta}), \vec{\zeta}]^2 d \vec{\zeta}} \overline{\Psi (\vec{x})} [\Psi (\vec{\eta}) \pm \Psi (-\vec{\eta})] d \vec{x} d \vec{\eta}\\
\end{eqnarray}
If then $\Psi$ is even, by choosing the negative sign in this equation, one obtains $\alpha - \overline{\alpha}=0$, whereas if $\Psi$ is odd, by choosing the plus sign, one finds $\alpha + \overline{\alpha}=0$. The eigenvalues of (\ref{eq:ii2}) associated with even eigenfunctions are real; the eigenvalues of (\ref{eq:ii2}) associated with odd eigenfunctions are pure imaginary. If follows then that (\ref{eq:ii2}) is equivalent to the pair of equations
\begin{eqnarray}
\beta_e \Psi_e(\vec{x}) = \int_{|\vec{x}| \leq a} Re(G_0 (0, \vec{x}; z, \vec{\eta}) e^{4 \delta^2 \sigma^2 k_0^4 \int K[(0, \vec{x}), (z, \vec{\eta}), \vec{\zeta}]^2 d \vec{\zeta}}) \Psi_e(\vec{\eta}) d \eta\\
\beta_o \Psi_o(\vec{x}) = \int_{|\vec{x}| \leq a} Im(G_0 (0, \vec{x}; z, \vec{\eta}) e^{4 \delta^2 \sigma^2 k_0^4 \int K[(0, \vec{x}), (z, \vec{\eta}), \vec{\zeta}]^2 d \vec{\zeta}}) \Psi_o(\vec{\eta}) d \eta
\end{eqnarray}
in which $\beta_e$ and $\beta_o$ are real. These equation have real symmetric kernels and we can fall back on the extensive theory in the literature treating such equations. It follows then from [22] that the eigenfunctions of (\ref{eq:ii2}) can be chosen real, orthogonal and complete in the class of functions square-integrable in the input domain. By iterating (\ref{eq:ii2}), one finds that the $\Psi$ also satisfy
\begin{eqnarray}
\nonumber &\ & |\alpha|^2 \Psi (\vec{x}) = \alpha (\overline{\alpha \Psi (\vec{x})})\\
\nonumber &=& \int_{|\vec{\eta}| \leq a} G_0 (0, \vec{x}; z, \vec{\eta}) e^{4 \delta^2 \sigma^2 k_0^4 \int K[(0, \vec{x}), (z, \vec{\eta}), \vec{\zeta}]^2 d \vec{\zeta}} \overline{\alpha \Psi (\vec{\eta})} d \vec{\eta}\\
\nonumber &=&\int_{|\vec{\eta}| \leq a} G_0 (0, \vec{x}; z, \vec{\eta}) e^{4 \delta^2 \sigma^2 k_0^4 \int K[(0, \vec{x}), (z, \vec{\eta}), \vec{\zeta}]^2 d \vec{\zeta}}\\
\nonumber &\ & \cdot \int_{|\vec{y}| \leq a} \overline{G_0 (0, \vec{\eta}; z, \vec{y})} e^{4 \delta^2 \sigma^2 k_0^4 \int \overline{K[(0, \vec{\eta}), (z, \vec{y}), \vec{\zeta}]}^2 d \vec{\zeta}} \overline{\Psi (\vec{y})} d \vec{y} d \vec{\eta}\\
&=& \int_{|\vec{y}| \leq a} K_r (\vec{x}, \vec{y}) \Psi (\vec{y}) d \vec{y} \label{eq:ii6}\\
&\ & \lambda ' = |\alpha|^2
\end{eqnarray}
where
\begin{equation}
K_r (\vec{x}, \vec{y}) = \int_{|\vec{\eta}| \leq a} G_0 (0, \vec{x}; z, \vec{\eta}) \overline{G_0 (0, \vec{\eta}; z, \vec{y})} e^{4 \delta^2 \sigma^2 k_0^4 \int (K[(0, \vec{x}), (z, \vec{\eta}), \vec{\zeta}]^2 + \overline{K[(0, \vec{\eta}), (z, \vec{y}), \vec{\zeta}]}^2) d \vec{\zeta}} d \vec{\eta}
\end{equation}
which is (\ref{eq:ii4}) in slightly altered notation when $a=b$, since $G_0 (0, \vec{\eta}; z, \vec{y}) = G_0 (0, \vec{y}; z, \vec{\eta})$. When $a \neq b$, it's a scaled version of (\ref{eq:ii4}) by paraxial approximation. Since the solution of (\ref{eq:ii2}) are complete, it follows that they are also a complete set of solutions of (\ref{eq:ii6}). As we asserted, to solve (\ref{eq:ii6}), it suffices to solve (\ref{eq:ii2}).

\newpage
\pagestyle{myheadings} 
\markright{  \rm \normalsize CHAPTER 9. \hspace{0.5cm} 
  Conclusion and Discussion}

\chapter{Conclusion and Discussion}
\thispagestyle{myheadings}

In the first part of this thesis, we have been dealing with the 3-D Helmholtz equation which is fundamental for all propagation theory. The Green's function for free-space propagation could be easily solved with specified boundary condition (Sommerfeld radiation condition).\\

In a turbulent medium where the magnitude of the air-density inhomogeneity is small, Rytov's method could be applied to approximate the Helmholtz equation for the frozen atmosphere. Within a limited but useful region of validity, both the solution and the Green's function of the Helmholtz equation are obtained explicitly.\\

Furthermore, by using the parabolic approximation, a simpler form could be obtained for both the perturbed wave field and the Green's function. In two specialized cases for the initial condition which are the plane wave case and beam wave case, the solution could be expressed in a form that numerical computation is easier to conducted.\\

In the second part of this thesis, through direct application of the extended Huygens-Fresnel principle, a general expression has been derived for the wave field received in the object plane for two optical problems - time reversal and apodization problem. In order to evaluate the ensemble average of the field and the intensity, we assume the refractive index disturbance to be white noise and take the Gaussian expectation on the inhomogeneous term. By doing so we end up with superposition integrals which can be computed by using numerical integration techniques.\\

At the very end of Chapter 8, we obtain an energy ratio (\ref{eq:ene}) which is very similar to the one for the free-space apodization problem except for an exponent term involved with the air-density inhomogeneity. The apodization problem for free-space propagation, which means $\delta=0$ in our case, has been well studied $^{[14][15][16]}$, and the free-space input and output eigenfunctions are prolate spheroidal wavefunctions. Therefore an immediate research problem based on this thesis is to build a connection between our results through Rytov's method and previous analysis of the analytic solution of the apodization problems.\\

To obtain the numerical results for these optical applications, we finally need to do integrations in 2-D or 3-D domains. In doing so, the parabolic approximation is helpful since it can separate the first dimension with the other two dimensions and thus reduce the complexity of the integrals. We can also specify the initial conditions in particular cases such as the plane wave and the beam wave.

\newpage
\pagestyle{myheadings} 
\markright{  \rm \normalsize BIBLIOGRAPHY. \hspace{0.5cm}
}

\thispagestyle{myheadings}

\end{spacing}

\end {document}